\documentclass[letterpaper, 10 pt, conference, onecolumn]{ieeeconf}  

\IEEEoverridecommandlockouts                              
\usepackage{tikz}

\usepackage{amsthm}    
\usepackage{amsmath} 
\usepackage{amssymb}  
\usepackage{amsfonts}
\usepackage{mathtools}
\usepackage{algpseudocode}
\usepackage{algorithm}
\usepackage{pgfplots}
\usepackage{multirow}
\usepackage{tcolorbox}
\usepackage{subcaption}
\usepackage{booktabs}

\pgfplotsset{compat=newest}
\usepgfplotslibrary{groupplots}

\usetikzlibrary{intersections}
\usetikzlibrary{calc}
\usetikzlibrary{positioning}
\usetikzlibrary{patterns, angles, quotes}
\usetikzlibrary{babel}
\usetikzlibrary{fadings,decorations.pathreplacing}
\usetikzlibrary{arrows.meta}

\usetikzlibrary{external}
\tikzset{external/only named=true}

\DeclareMathOperator{\diag}{diag}
\newtheorem{remark}{Remark}
\newtheorem{theorem}{Theorem}
\newtheorem{corollary}{Corollary}
\newtheorem{definition}{Definition}
\newtheorem{lemma}{Lemma}
\newtheorem{example}{Example}
\newtheorem{problem}{Problem}

\newcommand{\bbC}{\mathbb{C}}

\newcommand{\bbN}{\mathbb{N}}

\newcommand{\bbR}{\mathbb{R}}

\newcommand{\bA}{\boldsymbol{A}}
\newcommand{\bB}{\boldsymbol{B}}
\newcommand{\bC}{\boldsymbol{C}}
\newcommand{\bD}{\boldsymbol{D}}

\newcommand{\bI}{\boldsymbol{I}}

\newcommand{\bK}{\boldsymbol{K}}

\newcommand{\bP}{\boldsymbol{P}}
\newcommand{\bQ}{\boldsymbol{Q}}
\newcommand{\bR}{\boldsymbol{R}}
\newcommand{\bS}{\boldsymbol{S}}

\newcommand{\bW}{\boldsymbol{W}}

\newcommand{\calA}{\mathcal{A}}
\newcommand{\calB}{\mathcal{B}}
\newcommand{\calC}{\mathcal{C}}
\newcommand{\calD}{\mathcal{D}}
\newcommand{\calE}{\mathcal{E}}
\newcommand{\calF}{\mathcal{F}}

\newcommand{\calL}{\mathcal{L}}

\newcommand{\calS}{\mathcal{S}}

\newcommand{\calZ}{\mathcal{Z}}

\DeclareMathOperator{\CNN}{CNN}
\DeclareMathOperator{\FNN}{FNN}
\DeclareMathOperator{\HNN}{HNN}
\DeclareMathOperator*{\minimize}{minimize}

\newcommand{\lsb}{(}
\newcommand{\rsb}{)}

\title{\LARGE \bf
Convolutional Neural Networks as 2-D systems}

\author{Dennis Gramlich, Patricia Pauli, Carsten W. Scherer, Frank Allgöwer and Christian Ebenbauer
\thanks{The author thanks the International Max Planck Research
School for Intelligent Systems (IMPRS-IS) for support.}
\thanks{Dennis Gramlich and Christian Ebenbauer are with the Chair of Intelligent Control Systems,
	RWTH Aachen university,
	52074 Aachen, Germany
	{\tt\small \{dennis.gramlich,christian.ebenbauer\} @ic.rwth-aachen.de}}%
\thanks{Patricia Pauli and Frank Allgöwer are with the Institute for Systems Theory and Automatic Control, University of Stuttgart, 70569 Stuttgart, Germany
{\tt\small patricia.pauli@ist.uni-stuttgart.de}}%
\thanks{Carsten W. Scherer is with the Chair of Mathematical Systems Theory, University of Stuttgart, 70569 Stuttgart, Germany}
}

\begin{document}

\maketitle
\thispagestyle{empty}
\pagestyle{empty}

\begin{abstract}

This paper introduces a novel representation of convolutional Neural Networks (CNNs) in terms of 2-D dynamical systems. To this end, the usual description of convolutional layers with convolution kernels, i.e., the impulse responses of linear filters, is realized in state space as a linear time-invariant 2-D system. The overall convolutional Neural Network composed of convolutional layers and nonlinear activation functions is then viewed as 2-D version of a Lur´e system, i.e., a linear dynamical system interconnected with static nonlinear components. One benefit of this 2-D Lur´e system perspective on CNNs is that we can use robust control theory much more efficiently for Lipschitz constant estimation than previously possible.

\end{abstract}

\section{Introduction}


The success of Deep Learning is largely attributed to the evolution of Neural Network architectures. Structures like fully connected Neural Networks (FNN), recurrent Neural Networks (RNN), and convolutional Neural Networks have been the enablers of celebrated breakthroughs in science and are well accepted in, e.g., image and audio processing.


These breakthroughs are commonly attributed to properties such as invariances \cite{giles1987learning} and \emph{inductive biases} \cite{bietti2019inductive} of the Neural Network architectures. Inductive biases can be thought of as the implicit assumptions that come with the choice of a specific machine-learning model. An important example in the Neural Networks context is the inductive bias of \emph{deep architectures}, i.e., that a complicated nonlinear function is often a composition of multiple simpler ones \cite{goyal2022inductive}. 

For the FNN, the RNN, and the CNN, the deep architectures take the form of a chain of affine linear mappings and scalar, repeated nonlinearities (activation functions). The differences in these networks stem from the classes of affine linear mappings, and the vector spaces they operate on. Namely, in the FNN case, the affine linear mappings are not restricted and operate on finite dimensional vector spaces. In the RNN case, these mappings operate on a vector space of one-dimensional signals and are required to be invariant with respect to the shift operation along the time dimension of the signal. Finally, in the CNN case, the mappings operate (often) on a vector space of two-dimensional signals and are required to be invariant with respect to shifts along both space dimensions of the two-dimensional signal.

The latter two classes of linear mappings, i.e., linear shift-invariant mappings on signal spaces, are very well-researched in the domain of control theory, where we call them linear time-invariant (LTI) systems. It is thus not surprising that a lot of (ongoing) research has emerged linking control theory tools to Neural Networks. For RNNs, for example, long-term dependencies and the long-term behavior of the hidden state signals is an important aspect, and therefore, stability analysis from control theory is applied to address the well-known vanishing and exploding gradient problems \cite{pascanu2013difficulty}. Another interesting example is the analysis of Neural Networks based on robust control. It has been discovered in \cite{suykens1995artificial,levin1993control} that not just the single, e.g. recurrent or convolutional layer falls into the class of linear systems, but the whole Neural Network including the activation functions falls into the well-known class of Lur´e systems \cite{lur1944theory,gonzaga2012stability}. Moreover, the activation functions typically used in Neural Networks are well-studied nonlinearities in the context of Lur´e systems, namely sector-bounded and slope-restricted nonlinearities. This insight, which was first published in \cite{suykens1995artificial}, has been taken up again nearly three decades later in \cite{fazlyab2019efficient,fazlyab2020safety,hashemi2021certifying}, where the representation of Neural Networks as Lur´e systems is utilized to provide upper estimates on the Lipschitz constant of Neural Networks using robust control theory. The Lipschitz constant is an important property of Neural Networks, as it can be viewed as a measure of robustness for the Neural Network, i.e., the resilience of the Neural Network to adversarial attacks \cite{goodfellow2014explaining,szegedy2013intriguing}. Lipschitz constant estimation for Neural Networks and the training of Neural Networks with constraints on the Lipschitz constant are also investigated outside of the control engineering literature. For example, the result from \cite{anil2019sorting,tanielian2021approximating} is to be highlighted, where it is demonstrated that with the help of orthogonal layers and group-sort activation functions the approximation of any Lipschitz continuous function is possible. Of particular interest in this context are CNNs (see \cite{wang2020orthogonal}, where adversarial examples have been discovered for CNNs).

The Lur´e system perspective on Neural Networks has evolved to the utilization of dynamic integral quadratic constraints (IQCs) \cite{pauli2021linear,yin2021stability} for the stability analysis of neural controllers and other control application scenarios \cite{pauli2021offset}, and constrained FNN training \cite{pauli2021training,yin2021imitation,revay2021recurrent,revay2020lipschitz}. Recently, also the structure of RNNs as dynamical systems is being exploited in this direction of robust control analysis of Neural Networks \cite{fazlyab2019efficient} and, very recently, also the structure of 1-D CNNs \cite{pauli2022lipschitz}. 

What is not yet explored (to the best knowledge of the authors) is the possibility of analyzing Neural Networks with tools from the theory of 2-D dynamical systems. The paper \cite{pauli2022lipschitz} analyzes one-dimensional CNNs with tools from 1-D systems theory. The advantage of considering 2-D systems theory is, that, unlike \cite{pauli2022lipschitz}, which is restricted to CNNs for audio processing, the theory developed in this work is also applicable to CNNs for image processing. To this end, we view the layers of CNNs as 2-D systems. Indeed, the affine linear mappings in convolutional layers are finite impulse response (FIR) filters represented by their impulse response, i.e., their convolution kernel. Such filters can be viewed as n-D systems.

The domain of n-D systems is an extension of systems theory to signals in more than one independent variable, i.e., signals in multiple time and/or space dimensions. Such systems have been studied in the frequency domain for a long time (e.g. for filter design \cite{lin1996theory} or the analysis/discretization of PDEs \cite{willems2007stability}). The most well-known state space realizations, which are rather common in the control literature, have been proposed in 1975 \cite{roesser1975discrete} by Roesser and in 1980 \cite{fornasini1980stability} by Fornasini and Marchesini. In the present paper, we show that there is a 2-D systems state space realization for convolutional layers. As a consequence of this, we show in Section \ref{sec:2} that a (fully) convolutional Neural Network corresponds to the 2-D version of a Lur´e system. In Section \ref{sec:3}, we show how dissipativity theory and IQCs can be applied to 2-D system realizations of convolutional Neural Networks and use these findings to derive an effective linear matrix inequality problem formulation for estimating the Lipschitz constant of convolutional Neural Networks (Section \ref{sec:4}).

In particular, the proposed convex optimization problem for estimating the Lipschitz constant of a CNN should be more scalable than existing approaches, which require the embedding of CNNs into very high dimensional classes of feedforward Neural Networks with less structure. To achieve this convex formulation, IQCs are integrated into the setting of 2-D systems and 2-D Lur´e systems are introduced in this paper.




\section{Fully Convolutional Neural Networks and 2-D systems}
\label{sec:2}

Consider the space $\ell_{2e}^c(\bbN\times \bbN)$ of all signals $u : \bbN \times \bbN \to \bbR^c$, where $\bbN$ is defined as the set of all integers larger than or equal to zero. The signal space $\ell_2^c(\bbN \times \bbN)$ is defined as the subspace of $\ell_{2e}^c(\bbN\times \bbN)$ containing all signals with bounded energy, i.e., all functions $u \in \ell_{2e}^c(\bbN\times \bbN)$ with $\|u\|_{\ell_2}^2 = \sum_{i_1,i_2 = r}^\infty \|u\lsb i_1,i_2 \rsb \|^2 < \infty$. Here, $c \in \bbN$ denotes the number of channels (components) of the signal $u$. 
We view a (fully) convolutional Neural Network as an operator $\CNN_\theta : \ell_{2e}^{c_0}(\bbN \times \bbN) \to \ell_{2e}^{c_l}(\bbN \times \bbN)$ for some $c_0, c_l \in \bbN$ mapping infinite length two-dimensional signals to infinite length two-dimensional signals. At first glance, such a definition might seem odd, since convolutional Neural Networks usually take images, i.e., vectors from the finite-dimensional space $\bbR^{d_1 \times d_2 \times 3}$ as their input. However, this can be fully captured with our definition by embedding $\bbR^{d_1 \times d_2 \times 3}$ into the signal space $\ell_{2e}^3(\bbN \times \bbN)$ via
\begin{align}
    \calE: & \bbR^{d_1\times d_2 \times 3} \hookrightarrow \ell_{2e}^3(\bbN \times \bbN), v \mapsto (u\lsb i_1,i_2 \rsb )_{i_1,i_2 \in \bbN},\label{eq:image_embedding}
\end{align}
where the $i_3$-th entry $u_{i_3}\lsb i_1,i_2 \rsb$ of the vector $u\lsb i_1,i_2 \rsb$ is defined as
\begin{align*}
    u_{i_3}\lsb i_1,i_2 \rsb  = \begin{cases}
        v_{i_1,i_2,i_3} & i_1 \in \{1,\ldots,d_1\} \text{ and } i_2 \in \{1,\ldots,d_2\}\\
        0 & \text{else}
    \end{cases} 
\end{align*}
for $i_3 = 1,2,3$. By virtue of this embedding, our definition views images as a special case of 2-D signals also accounting for the fact that a CNN can process images of arbitrary size. We say that such a CNN has depth $l \in \bbN$ if the CNN is the concatenation of exactly $l$ convolutional layers and $l-1$ activation function layers. In this concatenation, the $k$-th convolutional layer can be defined as follows.

\begin{tcolorbox}[colback=gray!5!white,colframe=gray!75!black,title=Convolutional Layer]
A convolutional layer is described by a convolution kernel $\bK = (\bK\lsb i_1,i_2 \rsb )_{i_1,i_2 = 0}^{r}$ with $\bK\lsb i_1,i_2 \rsb \in \bbR^{c_k\times c_{k-1}}$ for $i_1,i_2 \in \{ 0,\ldots ,r\}$ and a bias $b \in \bbR^{c_k}$ as the mapping $\calC : \ell_{2e}^{c_{k-1}}(\bbN\times \bbN) \to \ell_{2e}^{c_k}(\bbN \times \bbN), (u\lsb i_1,i_2 \rsb ) \mapsto (y\lsb i_1,i_2 \rsb )$ defined by 
\begin{align}
    y\lsb i_1,i_2 \rsb  = b + \sum_{j_1,j_2 = 0}^{r} \bK\lsb j_1,j_2 \rsb  u\lsb i_1 - j_1,i_2 - j_2 \rsb \label{eq:convLayer}
\end{align}
for $i_1, i_2 \in \bbN$. The number $r$ defines the size of the convolution kernel. Mathematically, we define a convolutional layer in terms of the tuple of parameters $(\bK,b)$. Note that the convolution kernel $\bK$ implicitly defines the number of input channels $c_{k-1}$ and the number of output channels $c_k$.
\end{tcolorbox}
Clearly, the convolution kernel $\bK$ corresponds to the impulse response of the affine linear time-invariant filter \eqref{eq:convLayer}. We mention that it is not a restriction that we consider convolution kernels $\bK = (\bK\lsb i_1,i_2 \rsb )_{i_1,i_2 = r_-}^{r}$ with $r_- = 0$. If we are given a convolution kernel with $r_- < 0$, which is sometimes the case in the machine learning literature, then we can reindex this kernel and the output signal, such that $r_- = 0$ holds true. This is similar to software implementations, where the entries of $\bK$ are stored in arrays and negative array indices are uncommon.
In a CNN, convolutional layers, as above, are alternately concatenated with scalar non-linear functions $\sigma : \bbR \to \bbR$. To this end, the activation functions are identified with functions lifted to the signal space $\ell_{2e}^{c_k}(\bbN \times \bbN)$ which are given by element-wise application of $\sigma$, i.e.,
\begin{align*}
    \sigma : \ell_{2e}^{c_k}(\bbN \times \bbN) \to \ell_{2e}^{c_k}(\bbN\times \bbN), (y\lsb i_1,i_2 \rsb ) \mapsto (\sigma (y\lsb i_1,i_2 \rsb )).
\end{align*}
With these definitions of convolutional layers and activation functions, we are now able to provide a more detailed specification of CNNs.

\begin{tcolorbox}[colback=gray!5!white,colframe=gray!75!black,title=Fully convolutional Neural Network]
A convolutional Neural Network is a parameter dependent mapping $\CNN_\theta : \ell_{2e}^{c_0}(\bbN\times \bbN ) \to \ell_{2e}^{c_l}(\bbN\times \bbN)$ defined by
\begin{align}
    \CNN_\theta = \calC_{l}\circ \sigma \circ \calC_{l - 1} \circ \sigma \circ \cdots \circ \sigma \circ \calC_1, \label{eq:CNN}
\end{align}
where $\calC_1,\ldots,\calC_{l}$ are convolutional layers of compatible dimensions and $\sigma$ is a scalar activation function. The parameter $\theta$ contains the convolution kernels and biases $(b_1,\ldots,b_{l},\bK_1,\ldots,\bK_l)$ of all convolutional layers. Note that one convolutional Neural Network can process finite-length (finitely supported) input signals of different lengths.
\end{tcolorbox}
If we keep track of all the intermediate signals, we can also describe the convolutional Neural Network $\CNN_\theta : (u^{(0)}\lsb i_1,i_2 \rsb ) \mapsto (y^{(l)}\lsb i_1,i_2 \rsb)$ with the set of recursive equations
\begin{align*}
    (y^{(k)}\lsb i_1,i_2 \rsb ) &= \calC_k(u^{(k-1)}\lsb i_1,i_2 \rsb ) & k &=1,\ldots,l\\
    (u^{(k)}\lsb i_1,i_2 \rsb ) &= \sigma (y^{(k)}\lsb i_1,i_2 \rsb ) & k &=1,\ldots , l-1.
\end{align*}
We stress that many CNNs found in the literature contain additional types of layers such as batch/layer norm \cite{ba2016layer} or pooling layers \cite{gholamalinezhad2020pooling}, which are not captured by the definition in \eqref{eq:CNN}. It is for this reason that we refer to \eqref{eq:CNN} as \emph{fully convolutional neural network} \cite{long2015fully}. One of the most important additional layer types, namely, the (fully connected) linear layer, is captured in Section \ref{sec:4}. These linear layers are often appended at the end of the concatenation chain of a Neural Network when a 2-D signal (e.g. an image) shall be mapped to a \emph{non-signal} vector (such as class scores in image classification). Note that it should be possible to extend the present work to include pooling layers, as shown for example in \cite{pauli2022lipschitz} for 1-D convolutions.

The aim of this paper is twofold. The first goal is to show how a CNN can be modeled as a 2-D system (in state space). Existing work \cite{suykens1995artificial} already shows that ordinary fully connected Neural Networks can be represented as static Lur´e systems from Robust Control. However, these results do not yet generalize in a satisfying way to CNNs, where the structure of convolutions should be exploited to generate dynamic (Lur´e system) representations.

\begin{problem}
    \label{problem1}
    Derive a 2-D systems state space representation for convolutional Neural Networks.
\end{problem}

The second goal of the present work is to show one benefit of the 2-D state space representation by deriving tight bounds on the Lipschitz constant of CNNs.

\begin{problem}
    \label{problem2}
    For a given convolutional Network $\CNN_\theta$ with fixed parameter $\theta$, find an upper bound on the Lipschitz constant, i.e., find a value $\gamma \in \bbR$ such that
    \begin{align*}
        \| \CNN_\theta (u) - \CNN_\theta (u')\|_{\ell_2} &\leq \gamma \| u - u' \|_{\ell_2} & \forall u,u' \in \ell_2^{c_0}(\bbN \times \bbN).
    \end{align*}
\end{problem}

We emphasize that Problem \ref{problem2} aims at establishing an upper bound on the Lipschitz constant which applies for all 2-D signals from $\ell_2^{c_0}(\bbN\times \bbN)$. Particularly, $\gamma$ from Problem \ref{problem2} bounds the Lipschitz constant of the considered CNN for input images of \emph{any} size. This is a problem, which is only partly solved in the current literature. For a fixed input size, a CNN can be transformed into a fully connected Neural Network \cite{goodfellow2016deep} and methods for estimating the Lipschitz constant of fully connected Neural Networks exist (see \cite{fazlyab2019efficient}, \cite{fazlyab2020safety}, \cite{pauli2021training}). However, the representation of a CNN as a fully connected Neural Network is very inefficient for large input signals (it corresponds to representing a dynamical system by a highly sparse and repetitive Toeplitz matrix) and, therefore, the existing tools for computing the Lipschitz constant might become computationally intractable, while the new approach proposed in this paper does not. This inefficiency of Toeplitz matrix approaches is also reported in classical control theory where state space approaches have emerged as superior. The drawbacks of Toeplitz matrix/FIR representations of dynamical systems have been observed, for example, for Youla parameters which are often approximated as FIR filters \cite{scherer1999mixed,scherer2007multi} (the same is true for the system level parametrization). Furthermore, one clear advantage of CNNs, which is also often highlighted in the literature (see \cite{richter2021input}), is that they can process images of arbitrary size. This ability is clearly not reflected by studying the Toeplitz matrix as in the existing tools for Lipschitz constant estimation of CNNs. Unlike prior work on the subject, we handle input images of arbitrary size by tailoring our framework to convolutional mappings, establishing the fact that a convolutional layer such as \eqref{eq:convLayer} has a compact 2-D state space realization.

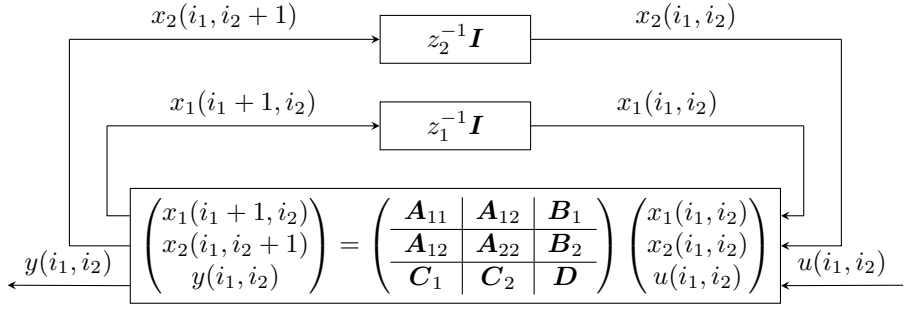
\begin{figure}
	\centering
	\begin{tikzpicture}
	\node [draw,
	minimum width=2cm,
	minimum height=1.2cm
	]  (2Dmatrix) at (0,0) {$\begin{pmatrix}
			x_1\lsb i_1+1,i_2 \rsb \\
			x_2\lsb i_1,i_2+1 \rsb \\
			y\lsb i_1,i_2 \rsb 
		\end{pmatrix}
		=
		\left(
		\begin{array}{c|c|c}
			\bA_{11} & \bA_{12} & \bB_1\\ \hline
			\bA_{12} & \bA_{22} & \bB_2\\ \hline
			\bC_1 & \bC_2 & \bD
		\end{array}
		\right)
		\begin{pmatrix}
			x_1\lsb i_1,i_2 \rsb \\
			x_2\lsb i_1,i_2 \rsb \\
			u\lsb i_1,i_2 \rsb 
		\end{pmatrix}
		$};
	
	\node [draw,
	minimum width=2cm,
	minimum height=0.6cm,
	above = 0.5cm of 2Dmatrix
	]  (firstshift) {$z_1^{-1} \bI$};
	
	\node [draw,
	minimum width=2cm,
	minimum height=0.6cm,
	above = 0.5cm of firstshift
	]  (secondshift) {$z_2^{-1} \bI$};
	
	\node[right = 3.5 cm of firstshift] (n1) {};
	\node[left = 3.5 cm of firstshift] (n2) {};
	\node[right = 4.0 cm of secondshift] (n3) {};
	\node[left = 4.0 cm of secondshift] (n4) {};
	\node[below right = 0.4cm and 1.5 cm of 2Dmatrix.east] (n5) {};
	\node[below left = 0.4cm and 1.5 cm of 2Dmatrix.west] (n6) {};
	
	\draw (firstshift.east) -- (n1.center) node[midway,above]{$x_1\lsb i_1,i_2 \rsb $};
	\draw[-stealth] (n2.center) -- (firstshift.west) node[midway,above]{$x_1\lsb i_1+1,i_2 \rsb $};
	\draw (secondshift.east) -- (n3.center) node[midway,above]{$x_2\lsb i_1,i_2 \rsb $};
	\draw[-stealth] (n4.center) -- (secondshift.west) node[midway,above]{$x_2\lsb i_1,i_2+1 \rsb $};
	
	\draw[-stealth] (n1.center) |- ($(2Dmatrix.east) + (0,0.4)$);
	\draw (n2.center) |- ($(2Dmatrix.west) + (0,0.4)$);
	\draw[-stealth] (n3.center) |- ($(2Dmatrix.east) + (0,0.0)$);
	\draw (n4.center) |- ($(2Dmatrix.west) + (0,-0.0)$);
	\draw[-stealth]  (n5.center) -- ($(2Dmatrix.east) - (0,0.53)$) node[midway,above]{$u\lsb i_1,i_2 \rsb $};
	\draw[-stealth]  ($(2Dmatrix.west) - (0,0.53)$) -- (n6.center) node[midway,above]{$y\lsb i_1,i_2 \rsb $};
	
\end{tikzpicture}
	\caption{Block diagram of a 2-D system as linear fractional representation in the time shift operators $z_1^{-1}$ and $z_2^{-1}$.}
	\label{fig:2DLFT}
\end{figure}

To obtain this state space realization, we apply the standard steps of first deriving the transfer function from the impulse response and then obtaining the state space realization from a linear fractional representation of the transfer function. Indeed, if we consider a two dimensional input signal $u\lsb i_1,i_2 \rsb $ and denote the shift operator along the index $i_1$ by $\sigma_1 : \ell_{2e} (\bbN \times \bbN) \to \ell_{2e}(\bbN \times \bbN), (u\lsb i_1,i_2 \rsb )_{i_1,i_2 \in \bbN} \mapsto (u\lsb i_1+1,i_2 \rsb )_{i_1,i_2 \in \bbN}$ and the shift operator along the index $i_2$ by $\sigma_2 : \ell_{2e} (\bbN \times \bbN) \to \ell_{2e}(\bbN \times \bbN), (u\lsb i_1,i_2 \rsb )_{i_1,i_2 \in \bbN} \mapsto (u\lsb i_1,i_2+1 \rsb )_{i_1,i_2 \in \bbN}$, then we can write the mapping $\calC$ defined in \eqref{eq:convLayer} as
\begin{align*}
	\calC : (u\lsb i_1,i_2 \rsb ) \mapsto \left(b + \left(\sum_{j_1,j_2 =0}^{r} \bK\lsb j_1,j_2 \rsb  \sigma_1^{-j_1} \sigma_2^{-j_2}\right) u\lsb i_1 ,i_2 \rsb \right).
\end{align*}
In this notation, $G(z_1,z_2) = \left(\sum_{j_1,j_2 = 0}^{r} \bK\lsb j_1,j_2 \rsb  z_1^{-j_1} z_2^{-j_2}\right)$ is equal to the transfer function of the filter defined by $\bK$. In this transfer function, $z_1^{-1}$ and $z_2^{-1}$ are complex variables, and, hence, $G(z_1,z_2)$ can be treated as a matrix of rational functions in $z_1^{-1}, z_2^{-1}$. According to \cite{bett1997linear}, such a transfer matrix can be realized in state space as a 2-D system. Indeed, viewing $G(z_1,z_2)$ as a matrix of rational functions, there always exists a linear fractional representation
\begin{align}
	G(z_1,z_2) = \begin{pmatrix}
		\bC_1 & \bC_2
	\end{pmatrix}
	\left(
	\begin{pmatrix}
		z_1 \bI & 0\\
		0 & z_2 \bI
	\end{pmatrix}
	-
	\begin{pmatrix}
		\bA_{11} & \bA_{12}\\
		\bA_{21} & \bA_{22}
	\end{pmatrix}
	\right)^{-1}
	\begin{pmatrix}
		\bB_1 \\ \bB_2
	\end{pmatrix}
	+ \bD \label{eq:tf_lfr}
\end{align}
composed of real matrices $(\bA_{11},\bA_{12},\bA_{21},\bA_{22}, \bB_1,\bB_2,\bC_1,\bC_2,\bD)$ of compatible dimensions, because $G$ has \emph{no poles in an environment of infinity}. It has been shown \cite{cockburn1997linear} that the problem of finding a realization for a 2-D system is a special case of finding a linear fractional representation (LFR) for a matrix depending rationally on complex variables. As it is common in Robust Control, we depict the linear fractional representation \eqref{eq:tf_lfr} as a block diagram in Figure~\ref{fig:2DLFT}. From this block diagram, we can directly see how to obtain a state space realization of this system as follows.

\begin{definition}[Roesser system]
	\label{def:RoesserSystem}
	An affine Roesser system is a 2-D system $\ell_2^{c_1}(\bbN\times \bbN) \to \ell_2^{c_2}(\bbN\times \bbN), (u\lsb i_1,i_2 \rsb ) \mapsto (y\lsb i_1,i_2 \rsb )$ with state space realization
	\begin{align}
		\begin{split}
			\begin{pmatrix}
				x_1\lsb i_1+1,i_2 \rsb \\
				x_2\lsb i_1,i_2+1 \rsb 
			\end{pmatrix}
			&=
			\begin{pmatrix}
				f_1\\
				f_2
			\end{pmatrix}
			+
			\begin{pmatrix}
				\bA_{11} & \bA_{12}\\
				\bA_{21} & \bA_{22}
			\end{pmatrix}
			\begin{pmatrix}
				x_1\lsb i_1,i_2 \rsb \\
				x_2\lsb i_1,i_2 \rsb 
			\end{pmatrix}
			+
			\begin{pmatrix}
				\bB_1\\
				\bB_2
			\end{pmatrix}u\lsb i_1,i_2 \rsb\\
			y\lsb i_1,i_2 \rsb &= g + \begin{pmatrix}
				\bC_1 & \bC_2
			\end{pmatrix}
			\begin{pmatrix}
				x_1\lsb i_1,i_2 \rsb \\
				x_2\lsb i_1,i_2 \rsb 
			\end{pmatrix}
			+ \bD u\lsb i_1,i_2 \rsb .
		\end{split}
		\label{eq:RoesserSys}
	\end{align}
	Here, the matrix tuple $(f_1,f_2,g,\bA_{11},\bA_{12},\bA_{21},\bA_{22},\bB_1,\bB_2,\bC_1,\bC_2, \bD)$ is called state space representation of the system, $x_1\lsb i_1,i_2 \rsb \in \bbR^{n_1}$, $x_2\lsb i_1,i_2 \rsb \in \bbR^{n_2}$ are the states, $u\lsb i_1,i_2 \rsb \in \bbR^{c_1}$ is the input and $y\lsb i_1,i_2 \rsb \in \bbR^{c_2}$ is the output of the system. We call \eqref{eq:RoesserSys} a linear time-invariant 2-D system if $f_1 = f_2 = 0$ and $g = 0$.
\end{definition}

The above state-space representation of 2-D systems has essentially been proposed by Roesser in \cite{roesser1975discrete}. Such state space realizations are not unique, which is shown e.g. by Fornasini and Marchesini, who propose an alternative and equally expressive realization in \cite{fornasini1980stability}. As we show in Example \ref{ex:2} the choice of the realization can make a difference for Lipschitz constant estimation.

A minor additions in Definition \ref{def:RoesserSystem} to the realization proposed by Roesser are the affine terms $f_1,f_2$ and $g$, which we introduce in order to deal with the bias terms in a Neural Network.

The arguments we put together up to this point essentially show the existence of a state space realization for convolutional layers, which is the statement of the following lemma.

\begin{lemma}[Convolutional Layers as 2-D systems]
    \label{lem:ConvLayer2Dsys}
    Let a convolution layer defined by a kernel $\bK = (\bK\lsb i_1,i_2 \rsb )_{i_1,i_2 = 0}^{r}$ and a bias $b$ be given. Then there exist matrices $f_1,f_2,g, \bA_{11},\bA_{12},\bA_{21},\bA_{22},\bB_1,\bB_2,\bC_1,\bC_2,\bD$ such that the convolutional layer \eqref{eq:convLayer} is realized by the 2-D system \eqref{eq:RoesserSys}. Here, realization means that the convolutional layer \eqref{eq:convLayer} and the 2-D system \eqref{eq:RoesserSys} define the same mapping from input signals $u \in \ell_{2e}^{c_1}(\bbN \times \bbN)$ to output signals $y \in \ell_{2e}^{c_2}(\bbN \times \bbN)$, when the states $x_1\lsb i_1,i_2 \rsb , x_2\lsb i_1,i_2 \rsb $ of \eqref{eq:RoesserSys} are initialized to $0$ for $(i_1,i_2) \in \{(i_1,0) \mid i_1 \in \bbN\}\cup \{(0,i_2) \mid i_2 \in \bbN\}$.
\end{lemma}

\begin{proof}
    Let $(u \lsb i_1,i_2 \rsb)$ be any signal from $\ell_2(\bbN \times \bbN)$. Then the $z$-transform $U(z_1,z_2) = \calZ ( u \lsb i_1,i_2 \rsb )$ converges for all $(z_1,z_2) \in \bbC$ with $|z_1|^2 + |z_2|^2 > 1$. The same is true for the output $(y_1 \lsb i_1,i_2 \rsb )$ of the convolutional layer, which satisfies $\calZ ( y_1 \lsb i_1,i_2 \rsb - b) = \widetilde{Y}_1(z_1,z_2) =  G(z_1,z_2) U(z_1,z_2)$, where the transfer function $G(z_1,z_2) = \left(\sum_{j_1,j_2 = 0}^{r} \bK\lsb j_1,j_2 \rsb  z_1^{-j_1} z_2^{-j_2}\right)$ of $\bK$ also converges in the considered region of convergence.

    Next, consider a 2-D system \eqref{eq:RoesserSys}, where the matrices $\bA_{11},\bA_{12},\bA_{21},\bA_{22},\bB_1,\bB_2,\bC_1,\bC_2,\bD$ are obtained from a linear fractional representation \eqref{eq:tf_lfr} of $G$ and the remaining vectors $f_1, f_2, g$ are set to $f_1 = f_2 = 0$ and $g = b$. If the states $x_1\lsb i_1,i_2 \rsb , x_2\lsb i_1,i_2 \rsb $ of this system are initialized to $0$ for $(i_1,i_2) \in \{(i_1,0) \mid i_1 \in \bbN\}\cup \{(0,i_2) \mid i_2 \in \bbN\}$ and $(u \lsb i_1,i_2 \rsb)$ is considered as input, then the $z$-transforms of the states $X_1(z_1,z_2) = \calZ (x_1\lsb i_1,i_2 \rsb )$, $X_2(z_1,z_2) = \calZ(x_2(i_1,i_2))$ converge for a sufficiently large radius $\rho\geq 1$ and $|z_1|^2 + |z_2|^2 > \rho^2$. Their $z$-transforms satisfy
    \begin{align*}
        \begin{pmatrix}
			z_1 \bI & 0\\
			0 & z_2 \bI
		\end{pmatrix}
				\begin{pmatrix}
					X_1(z_1,z_2)\\
					X_2(z_1,z_2)
				\end{pmatrix}
				&=
				\begin{pmatrix}
					\bA_{11} &\bA_{12}\\
					\bA_{21} & \bA_{22}
				\end{pmatrix}
				\begin{pmatrix}
					X_1(z_1,z_2)\\
					X_2(z_1,z_2)
				\end{pmatrix}
				+
				\begin{pmatrix}
					\bB_1\\
					\bB_2
				\end{pmatrix} U(z_1,z_2),
    \end{align*}
    which implies
    \begin{align*}
	\begin{pmatrix}
		X_1(z_1,z_2)\\
		X_2(z_1,z_2)
	\end{pmatrix}
	=
	\left(
	\begin{pmatrix}
		z_1 \bI & 0\\
		0 & z_2 \bI
	\end{pmatrix}
	-
	\begin{pmatrix}
		\bA_{11} & \bA_{12}\\
		\bA_{21} & \bA_{22}
	\end{pmatrix}
	\right)^{-1}
	\begin{pmatrix}
		\bB_1\\
		\bB_2
	\end{pmatrix} U(z_1,z_2).
    \end{align*}
    Next, we utilize the output equation of \eqref{eq:RoesserSys} to obtain that $\widetilde{Y}_2 = \calZ (y_2 \lsb z_1,z_2 \rsb - b)$ satisfies
    \begin{align}
		\widetilde{Y}_2(z_1,z_2) &= \begin{pmatrix}
			\bC_1 & \bC_2
		\end{pmatrix}
		\left(
		\begin{pmatrix}
			z_1 \bI & 0\\
			0 & z_2 \bI
		\end{pmatrix}
		-
		\begin{pmatrix}
			\bA_{11} & \bA_{12}\\
			\bA_{21} & \bA_{22}
		\end{pmatrix}
		\right)^{-1}
		\begin{pmatrix}
			\bB_1\\
			\bB_2
		\end{pmatrix} U(z_1,z_2)
		+ \bD U(z_1,z_2), \label{eq:linfrac2}
    \end{align}
    where $(y_2\lsb i_1,i_2\rsb )$ is the output of the 2-D system. Observing that \eqref{eq:linfrac2} contains the linear fractional representation of $G$, we can conclude that $\widetilde{Y}_2(z_1,z_2) = \widetilde{Y}(z_1,z_2)$ holds true for all $z_1,z_2 \in \bbC$ with $|z_1|^2 + |z_2|^2 > \rho^2$. By the identity theorem for holomorphic functions, this implies that $y_1 \lsb i_1,i_2\rsb - b$ and $y_2\lsb i_1,i_2\rsb - b$ are equal.
    
    To generalize the result to $u \in \ell_{2e}(\bbN \times \bbN)$, we can utilize standard truncation arguments (see \cite{scherer2000linear} for 1-D systems) which utilize the causality of the studied systems.
\end{proof}

With Lemma \ref{lem:ConvLayer2Dsys}, we establish that a convolutional layer can be represented as a 2-D system. This result could also be further utilized to estimate the Lipschitz constant of a single convolutional layer. To also handle CNNs, which consist of convolutional layers \emph{and} nonlinear activation functions, we extend the definition of affine Roesser systems by including nonlinearities. The resulting new system class can be thought of as a 2-D version of the Lur´e system.

\begin{definition}[Affine 2-D Lur´e systems]
    \label{def:LureSys}
    A 2-D Lur´e system is defined by real matrices $f_1$, $f_2$, $g_1$, $g_2$, $\calA_{11}$, $\calA_{12}$, $\calA_{21}$, $\calA_{22}$, $\calB_{11}$, $\calB_{12}$, $\calB_{21}$, $\calB_{22}$, $\calC_{11}$, $\calC_{12}$, $\calC_{21}$, $\calC_{22}$, $\calD_{11}$, $\calD_{12}$, $\calD_{21}$ and $\calD_{22}$ of compatible dimensions, and a nonlinear function $\phi : \bbR^d \to \bbR^d$ by the signal relations
    \begin{align}
        \begin{pmatrix}
        x_1\lsb i_1+1,i_2 \rsb \\
        x_2\lsb i_1,i_2+1 \rsb \\
        z\lsb i_1,i_2 \rsb \\
        y\lsb i_1,i_2 \rsb 
        \end{pmatrix}
        =
        \begin{pmatrix}
        	f_1 & \calA_{11} & \calA_{12} & \calB_{11} & \calB_{12} \\
	        f_2 & \calA_{21} & \calA_{22} & \calB_{21} & \calB_{22} \\
	        g_1 & \calC_{11} & \calC_{12} & \calD_{11} & \calD_{12}\\
	        g_2 & \calC_{21} & \calC_{22} & \calD_{21} & \calD_{22}
        \end{pmatrix}
        \begin{pmatrix}
        1\\
        x_1\lsb i_1,i_2 \rsb \\
        x_2\lsb i_1,i_2 \rsb \\
        w\lsb i_1,i_2 \rsb \\
        u\lsb i_1,i_2 \rsb
        \end{pmatrix}, \hspace{5mm} w\lsb i_1,i_2 \rsb  = \phi (z\lsb i_1,i_2 \rsb ). \label{eq:2D_linearFractionalRepresentation}
    \end{align}
    Here, $x_1\lsb i_1,i_2 \rsb \in \bbR^{n_1}$, $x_2\lsb i_1,i_2 \rsb \in \bbR^{n_2}$ are the states, $w\lsb i_1,i_2 \rsb , z\lsb i_1,i_2\rsb \in \bbR^d$ are the input and output of the nonlinearity $\phi$, $u\lsb i_1,i_2 \rsb \in \bbR^{c_1}$ is the input of the system and $y\lsb i_1,i_2 \rsb \in \bbR^{c_2}$ is the output of the system. This is a feedback interconnection between an affine LTI 2-D system and the nonlinear map $\phi$. We call the affine 2-D Lur´e system a (linear) 2-D Lur´e system if $f_1 = f_2 = 0$, $g_1 = 0$, and $g_2 = 0$.
\end{definition}

In these 2-D Lur´e systems, the nonlinearities $\phi$ can be included in a feedback loop with an affine Roesser system which turns out to be sufficient for describing the CNN \eqref{eq:CNN}. Like Definition \ref{def:RoesserSystem}, Definition \ref{def:LureSys} includes the affine terms $f_1$, $f_2$, $g_1$ and $g_2$ to account for the biases in a Neural Network.
To show that a CNN \eqref{eq:CNN} can be expressed in the form \eqref{eq:2D_linearFractionalRepresentation}, we utilize that every layer $\calC_i$, $i = 1,\ldots,l$ has a 2-D state space representation $(f_1^{(i)},f_2^{(i)},g^{(i)},\bA_{11}^{(i)},\bA_{12}^{(i)},\bA_{21}^{(i)},\bA_{22}^{(i)},\bB_1^{(i)},\bB_2^{(i)},\bC_{1}^{(i)},\bC_2^{(i)},\bD^{(i)})$ (Lemma \ref{lem:ConvLayer2Dsys}). With the aid of these matrices, we can generate a Lur´e system \eqref{eq:2D_linearFractionalRepresentation} by simply writing down the equations of all layers in the vectorized form. We obtain
\begin{align}
    \overbrace{
    \begin{pmatrix}
    x_1^{(1)}\lsb i_1+1,i_2 \rsb \\
    x_1^{(2)}\lsb i_1+1,i_2 \rsb \\
    \vdots\\
    x_1^{(l)}\lsb i_1+1,i_2 \rsb 
    \end{pmatrix}}^{x_1\lsb i_1+1,i_2 \rsb }
    &=
    \overbrace{
    \begin{pmatrix}
    \bA_{11}^{(1)} & & &\\
    & \bA_{11}^{(2)} & &\\
    & & \ddots & \\
    & & & \bA_{11}^{(l)}
    \end{pmatrix}}^{\calA_{11}}
    \overbrace{
    \begin{pmatrix}
    x_1^{(1)}\lsb i_1,i_2 \rsb \\
    x_1^{(2)}\lsb i_1,i_2 \rsb \\
    \vdots\\
    x_1^{(l)}\lsb i_1,i_2 \rsb 
    \end{pmatrix}}^{x_1\lsb i_1,i_2 \rsb }
    +
    \overbrace{
    \begin{pmatrix}
    \bA_{12}^{(1)} & & &\\
    & \bA_{12}^{(2)} & &\\
    & & \ddots & \\
    & & & \bA_{12}^{(l)}
    \end{pmatrix}}^{\calA_{12}}
    \overbrace{
    \begin{pmatrix}
    x_2^{(1)}\lsb i_1,i_2 \rsb \\
    x_2^{(2)}\lsb i_1,i_2 \rsb \\
    \vdots\\
    x_2^{(l)}\lsb i_1,i_2 \rsb 
    \end{pmatrix}}^{x_2\lsb i_1,i_2-1 \rsb }\nonumber\\
    &+
    \underbrace{
    \begin{pmatrix}
    0 & \cdots & 0\\
    \bB_1^{(2)} & &\\
    & \ddots &\\
    & & \bB_1^{(l)}
    \end{pmatrix}}_{\calB_{11}}
    \underbrace{
    \begin{pmatrix}
    u^{(2)}\lsb i_1,i_2 \rsb \\
    \vdots\\
    u^{(l)}\lsb i_1,i_2 \rsb 
    \end{pmatrix}}_{w\lsb i_1,i_2 \rsb }
    +
    \underbrace{
    \begin{pmatrix}
    \bB_1^{(1)}\\
    0\\
    \vdots\\
    0
    \end{pmatrix}}_{\calB_{12}}
    \underbrace{
    u^{(1)}\lsb i_1,i_2 \rsb }_{u\lsb i_1,i_2 \rsb } \label{eq:realizationCNN1}
\end{align}
for the states $x_1^{(i)}\lsb i_1+1,i_2 \rsb $, $i = 1,\ldots,l$ and a similar equation for $x_2^{(i)}\lsb i_1,i_2+1 \rsb $. We further obtain
\begin{align}
    \overbrace{
    \begin{pmatrix}
    y^{(1)}\lsb i_1,i_2 \rsb \\
    \vdots\\
    y^{(l-1)}\lsb i_1,i_2 \rsb 
    \end{pmatrix}}^{z\lsb i_1,i_2 \rsb }
    &=
    \overbrace{
    \begin{pmatrix}
    \bC_1^{(1)} & & &\\
    & \ddots & &\\
    & & \bC_1^{(l-1)} & 0
    \end{pmatrix}}^{\calC_{11}}
    \overbrace{
    	\begin{pmatrix}
    		x_1^{(1)}\lsb i_1,i_2 \rsb \\
    		\vdots\\
    		x_1^{(l-1)}\lsb i_1,i_2 \rsb \\
    		x_1^{(l)}\lsb i_1,i_2 \rsb 
    \end{pmatrix}}^{x_1\lsb i_1,i_2 \rsb }
    +
    \overbrace{
    \begin{pmatrix}
    \bC_2^{(1)} & & &\\
    & \ddots & &\\
    & & \bC_2^{(l-1)} & 0
    \end{pmatrix}}^{\calC_{12}}
    \overbrace{
    	\begin{pmatrix}
    		x_2^{(1)}\lsb i_1,i_2 \rsb \\
    		\vdots\\
    		x_2^{(l-1)}\lsb i_1,i_2 \rsb \\
    		x_2^{(l)}\lsb i_1,i_2 \rsb 
    \end{pmatrix}}^{x_2\lsb i_1,i_2 \rsb }\nonumber\\
    &+
    \underbrace{
    \begin{pmatrix}
    0 & & & \\
    \bD^{(2)} & 0 & \\
    & \ddots & \ddots & \\
    & & \bD^{(l-1)} & 0
    \end{pmatrix}}_{\calD_{11}}
    \underbrace{
    	\begin{pmatrix}
    		u^{(2)}\lsb i_1,i_2 \rsb \\
    		\vdots\\
    		u^{(l-1)}\lsb i_1,i_2 \rsb \\
    		u^{(l)}\lsb i_1,i_2 \rsb 
    \end{pmatrix}}_{w\lsb i_1,i_2 \rsb }
    +
    \underbrace{
    \begin{pmatrix}
    \bD^{(1)}\\
    0\\
    \vdots\\
    0
    \end{pmatrix}}_{\calD_{12}}
    \underbrace{
    	u^{(1)}\lsb i_1,i_2 \rsb }_{u\lsb i_1,i_2 \rsb }
    +
	\underbrace{
	\begin{pmatrix}
		g^{(1)}\\
		\vdots\\
		g^{(l-1)}
	\end{pmatrix}}_{g_1} \label{eq:realizationCNN2}
\end{align}
for the outputs $y^{(i)}\lsb i_1,i_2 \rsb $, $i = 1,\ldots, l-1$ and
\begin{align}
    \overbrace{
    y^{(l)}\lsb i_1,i_2 \rsb }^{y\lsb i_1,i_2 \rsb }
    &=
    \overbrace{
    \begin{pmatrix}
    0 & \cdots & 0 & \bC_1^{(l)}
    \end{pmatrix}}^{\calC_{21}}
    \overbrace{
    	\begin{pmatrix}
    		x_1^{(1)}\lsb i_1,i_2 \rsb \\
    		\vdots\\
    		x_1^{(l-1)}\lsb i_1,i_2 \rsb \\
    		x_1^{(l)}\lsb i_1,i_2 \rsb 
    \end{pmatrix}}^{x_1\lsb i_1,i_2 \rsb }
    +
    \overbrace{
    \begin{pmatrix}
    0 & \cdots & 0 & \bC_2^{(l)}
    \end{pmatrix}}^{\calC_{22}}
    \overbrace{
    	\begin{pmatrix}
    		x_2^{(1)}\lsb i_1,i_2 \rsb \\
    		\vdots\\
    		x_2^{(l-1)}\lsb i_1,i_2 \rsb \\
    		x_2^{(l)}\lsb i_1,i_2 \rsb 
    \end{pmatrix}}^{x_2\lsb i_1,i_2 \rsb } \nonumber\\
    &+
    \underbrace{
    \begin{pmatrix}
    0 & \cdots & 0 & \bD^{(l)}
    \end{pmatrix}}_{\calD_{21}}
    \underbrace{
    	\begin{pmatrix}
    		u^{(2)}\lsb i_1,i_2 \rsb \\
    		\vdots\\
    		u^{(l-1)}\lsb i_1,i_2 \rsb \\
    		u^{(l)}\lsb i_1,i_2 \rsb 
    \end{pmatrix}}_{w\lsb i_1,i_2 \rsb }
    +
    \underbrace{
    0}_{\calD_{22}} 
	\underbrace{
		u^{(1)}\lsb i_1,i_2 \rsb }_{u\lsb i_1,i_2 \rsb }
	+
	\underbrace{g^{(l)}}_{g_2} \label{eq:realizationCNN3}
\end{align}
for $y^{(l)}\lsb i_1,i_2 \rsb $. Finally, we can summarize the activation function layers in
\begin{align*}
	\begin{pmatrix}
		u^{(2)}\lsb i_1,i_2 \rsb \\
		\vdots\\
		u^{(l)}\lsb i_1,i_2 \rsb 
	\end{pmatrix}
	&=
	\sigma
	\begin{pmatrix}
		y^{(1)}\lsb i_1,i_2 \rsb \\
		\vdots\\
		y^{(l-1)}\lsb i_1,i_2 \rsb 
	\end{pmatrix} \Leftrightarrow w\lsb i_1,i_2 \rsb  = \sigma (z\lsb i_1,i_2 \rsb ).
\end{align*}

With the newly defined variables $f_1,f_2,g_1,g_2,\calA_{11},\ldots, \calD_{22}$, we can pose the overall convolutional Neural Network as a single 2-D system interconnected with the nonlinear function $\phi = \sigma$ as in \eqref{eq:2D_linearFractionalRepresentation}. We further highlight the structure of the obtained 2-D Lur´e system, which results from the series interconnection of the convolutional layers and diagonal nonlinearities (diagonal/sub-diagonal matrices in \eqref{eq:2D_linearFractionalRepresentation}). This structure enables us for any given input $u\lsb i_1,i_2 \rsb  = u^{(1)}\lsb i_1,i_2 \rsb $ to first compute the states and outputs of $\calC_1$, subsequently compute the input to $\calC_2$ by $u^{(2)}\lsb i_1,i_2 \rsb = \sigma (y^{(1)}\lsb i_1,i_2 \rsb )$, then compute states and outputs of $\calC_2$ and, going on in this fashion, finally compute a solution to the equation system \eqref{eq:2D_linearFractionalRepresentation}. This shows that the obtained 2-D Lur´e system (from the CNN \eqref{eq:CNN}) is well-posed, i.e., for any input $u\lsb i_1,i_2 \rsb  \in \ell_{2e}^{c_0}(\bbN\times \bbN)$ there exist state and output signals satisfying the state and output equations. We summarize our result in the following theorem.

\begin{theorem}[CNNs as affine 2-D Lur´e systems]
    \label{thm:CNNasLureSys}
    Let a convolutional Neural Network $\CNN_\theta = \calC_l \circ \sigma \circ \calC_{\ell - 1} \circ \cdots \circ \sigma \circ \calC_1$ be given, where $\theta = (b_1,\ldots,b_l,\bK_1,\ldots,\bK_l)$. Then there exist matrices $f_1,f_2,g_1,g_2,\calA_{11}, \ldots , \calD_{22}$, such that the CNN \eqref{eq:CNN} can be realized as a well-posed Lur´e system of the form \eqref{eq:2D_linearFractionalRepresentation} with the nonlinearity $\phi = \sigma$.
\end{theorem}

\begin{proof}
    The derivation above the theorem is a constructive proof of this result.
\end{proof}

We note that \emph{naive} 2-D system representations of given CNNs \eqref{eq:CNN} can be constructed intuitively. The following example shows this for a single convolutional layer. The realization for the total CNN can then be obtained from the equations \eqref{eq:realizationCNN1}, \eqref{eq:realizationCNN2}, \eqref{eq:realizationCNN3}.

\begin{example}[Realizing a convolutional layer]
    \label{ex:realization}
    Let a convolutional layer \eqref{eq:convLayer} represented by $(\bK,b)$ be given. To obtain a realization \eqref{eq:RoesserSys}, we set $g = b$ and $f_1 = f_2 = 0$. If the states $x_1\lsb i_1,i_2 \rsb $ and $x_2\lsb i_1,i_2 \rsb $ are defined as
    \begin{align*}
        x_1\lsb i_1,i_2 \rsb  &\bumpeq \begin{pmatrix}
            u\lsb i_1-1,i_2 \rsb  & \cdots & u\lsb i_1-r,i_2 \rsb \\
            \vdots & \ddots & \vdots\\
            u\lsb i_1-1,i_2-r \rsb  & \cdots & u\lsb i_1-r,i_2-r \rsb 
        \end{pmatrix}, &
        x_2\lsb i_1,i_2 \rsb  &\bumpeq \begin{pmatrix}
            u\lsb i_1,i_2-1 \rsb  & \cdots & u\lsb i_1,i_2-r \rsb \\
            \vdots & \ddots & \vdots\\
            u\lsb i_1-r,i_2-1 \rsb  & \cdots & u\lsb i_1-r,i_2-r \rsb 
        \end{pmatrix},
    \end{align*}
    where $\bumpeq$ denotes that the left side is the vectorization of the matrix on the right, then we can determine the output $y\lsb i_1,i_2 \rsb $ from the states and input. Namely, since these states $x_1\lsb i_1,i_2 \rsb , x_2\lsb i_1,i_2 \rsb $ contain all the $u$-terms in $\sum_{j_1,j_2 = 0}^{r} \bK\lsb j_1,j_2 \rsb  u\lsb i_1 - j_1,i_2 - j_2 \rsb $ except for $u\lsb i_1,i_2 \rsb $, we can easily choose $\bC_1$, $\bC_2$ and $\bD$, such that
    \begin{align*}
        y\lsb i_1,i_2 \rsb  = b + \sum_{j_1,j_2 = 0}^{r} \bK\lsb j_1,j_2 \rsb  u\lsb i_1 - j_1,i_2 - j_2 \rsb  = b + \bC_1 x_1\lsb i_1,i_2 \rsb  + \bC_2 x_2\lsb i_1,i_2 \rsb  + \bD u\lsb i_1,i_2 \rsb .
    \end{align*}
    Choosing $\bA_{11}$, $\bA_{12}$ and $\bB_1$ as \emph{selection} matrices that shift certain entries to new positions and zero to other positions
    \begin{align*}
        \bA_{11} x_1\lsb i_1,i_2 \rsb  &\bumpeq \begin{pmatrix}
            0 & u\lsb i_1-1,i_2 \rsb  & \cdots & u\lsb i_1+1-r,i_2 \rsb \\
            \vdots & \vdots & \ddots & \vdots\\
            0 & u\lsb i_1-1,i_2-r \rsb  & \cdots & u\lsb i_1+1-r,i_2-r \rsb 
        \end{pmatrix},\\
        \bA_{12} x_2\lsb i_1,i_2 \rsb  &\bumpeq
        \begin{pmatrix}
            0 & 0 & \cdots & 0\\
            u\lsb i_1,i_2-1 \rsb & 0 & \cdots & 0\\
            \vdots & \vdots & \ddots & \vdots\\
            u\lsb i_1,i_2-r \rsb  & 0 & \cdots & 0
        \end{pmatrix}, \qquad
        \bB_1 u\lsb i_1,i_2 \rsb  \bumpeq
        \begin{pmatrix}
            u\lsb i_1,i_2 \rsb  & 0 & \cdots & 0\\
            0 & 0 & \cdots & 0\\
            \vdots & \vdots & \ddots & \vdots\\
            0 & 0 & \cdots & 0
        \end{pmatrix},
    \end{align*}
    we can also satisfy the state equations in \eqref{eq:RoesserSys} with these states. Here, $\bA_{11}$ shifts the columns of the \emph{matrix representation} of $x_1$ to the right, $\bA_{12}$ moves the first row of $x_2$ into the first column of its output with a zero on top and all remaining columns equal to zero and $\bB_1$ simply moves its input to the $(1,1)$ position of its output. The exact same choices for the remaining matrices, i.e., $\bA_{22} = \bA_{11}$, $\bA_{21} = \bA_{12}$, $\bB_2 = \bB_1$ render the state equations in \eqref{eq:RoesserSys} satisfied. Intuitively, $\bA_{11}$ and $\bA_{22}$ should be viewed as shift operators along $i_1$ and $i_2$ and the matrices $\bA_{12}$ and $\bB_1$ (or $\bA_{21}$ and $\bB_2$ respectively) fill in the missing data.
\end{example}

It should be noted that the realization from Example \ref{ex:realization} introduces $c_1r(r+1)$ states for $x_1$ and the same number of states for $x_2$, which can be much larger than the minimal number of states required for the realization of \eqref{eq:convLayer}. Indeed, in \cite{varoufakis1987minimal} it is shown that for the case of an FIR filter, the number of states can be as small as $c_1r$ for both $x_1$ and $x_2$. Furthermore, the realization as exemplified in Example \ref{ex:realization} is not controllable in the sense of 2-D systems \cite{klamka2005controllability}, i.e., there exist states in this realization which can never be reached. We mention that there are many notions of controllability for 2-D systems \cite{eising1979controllability,rocha1991controllability} and that the equivalence between controllability and observability, and minimal realization does not hold for these notions, as it does in the 1-D system case. Moreover, unlike for 1-D systems, the problem of minimal realization of general transfer functions is unsolved for 2-D systems \cite{mentzelopoulou1991n,antoniou1986state,galkowski1997elementary,kung1977new}. Still, there are tools for generating \emph{small} realizations of 2-D systems, which realize FIR filters minimally \cite{varoufakis1987minimal}. For this reason, we recommend realizing a convolutional layer by first formulating the transfer function of \eqref{eq:convLayer} and then using e.g. the toolbox \cite{biannic2016generalized} to derive a smaller 2-D system realization of the individual layers than provided in Example \ref{ex:realization}. The overall Lur´e system should then again be constructed as in \eqref{eq:realizationCNN1}, \eqref{eq:realizationCNN2}, \eqref{eq:realizationCNN3}.

\begin{figure}
    \centering
    \newcommand*{\xMin}{0}%
\newcommand*{\xMax}{6}%
\newcommand*{\yMin}{0}%
\newcommand*{\yMax}{6}%
\newcommand*{\iOne}{3}%
\newcommand*{\iTwo}{4}%
\usetikzlibrary{math}
\begin{tikzpicture}[scale=0.80]
	\foreach \i in {\xMin,...,\xMax} {
		\draw [very thin,gray] (\i,-\yMin) -- (\i,-\yMax) ;
	}
	\foreach \i in {\yMin,...,\yMax} {
		\draw [very thin,gray] (\xMin,-\i) -- (\xMax,-\i);
	}
	
	\foreach \i in {\xMin,...,\xMax} {
		\foreach \j in {\yMin,...,\yMax} {
			\node[circle,inner sep=0.05cm,draw,minimum size=0.02cm,fill] (x$\i\j$) at ($(\i,-\j)$) {};
			\pgfmathparse{\i + 1 + 7*\j}
			\xdef\n{\pgfmathresult}
			\node [below right=0.0cm and -0.05cm of x$\i\j$] {\tiny $u(\i,\j)$};
		}
	}

	\draw[-stealth,thick]  (0,0) -- (7,0) node[at end,above]{$i_1$};
        \draw[-stealth,thick]  (0,0) -- (0,-7) node[at end,left]{$i_2$};

        \draw[green] (2,-4) -- (0,-4);
        \draw[green] (2,-4) -- (2,0);
        \fill [fill=green,fill opacity = 0.1] (0,0) rectangle (2,-4);

        \draw[blue] (3,-3) -- (0,-3);
        \draw[blue] (3,-3) -- (3,0);
        \fill [fill=blue,fill opacity = 0.1] (0,0) rectangle (3,-3);

        \draw[red] (3.05,-4.05) -- (0,-4.05);
        \draw[red] (3.05,-4.05) -- (3.05,0);
        \fill [fill=red,fill opacity = 0.1] (2,-3) rectangle (3,-4);
\end{tikzpicture}
    \caption{
    This figure shows the causal dependence of the states $x_1\lsb 2,4 \rsb$ and $x_2\lsb 3,3 \rsb$ on inputs $u\lsb j_1,j_2 \rsb$ in a 2-D dynamical system. The green cone represents the inputs $x_1\lsb 2,4 \rsb$ may depend on, while the blue cone illustrates the inputs that influence $x_2\lsb 3,3 \rsb$. The combination of the green and blue cones with the area for $u\lsb 3,4 \rsb$ results in the red cone, which encompasses all inputs that affect the states $x_1\lsb 3,4 \rsb$ and $x_2\lsb 3,4 \rsb$. The states $x_1\lsb 3,4 \rsb$ and $x_2\lsb 3,4 \rsb$ should be understood as a concise summary of the information contained in the red cone, which is provided to future states in the system.}
    \label{fig:causality}
\end{figure}
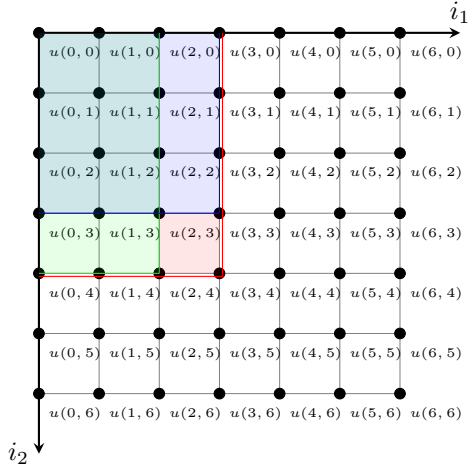

\begin{remark}
    We note that Example \ref{ex:realization} can also help in understanding the causality properties of 2-D systems. This causality requires that $x_1\lsb i_1,i_2 \rsb $ must depend only on the inputs $u\lsb j_1,j_2 \rsb $ for $j_1 < i_1$ and $j_2 \leq i_2$. Similarly, $x_2\lsb i_1,i_2 \rsb $ must depend only on $u\lsb j_1,j_2 \rsb $ for $j_1 \leq i_1$ and $j_2 < i_2$. The union of the \emph{causality cones} of $x_1\lsb i_1,i_2 \rsb $ and $x_2\lsb i_1,i_2 \rsb$ with $u\lsb i_1,i_2 \rsb $ yields the cone of all $(j_1,j_2)$ with $j_1 \leq i_1$ and $j_2 \leq i_2$, which corresponds exactly to the cone of inputs $x_1\lsb i_1+1,i_2 \rsb $ and $x_2\lsb i_1,i_2+1 \rsb $ may depend on. Moreover, this cone contains all the information to compute $y\lsb i_1,i_2 \rsb $. We have illustrated this reasoning in Figure \ref{fig:causality}. In the case of Example \ref{ex:realization}, we have simply chosen the states $x_1\lsb i_1,i_2 \rsb $ and $x_2\lsb i_1,i_2 \rsb $ as all relevant inputs in their respective causality cones, i.e., all $u\lsb j_1,j_2 \rsb $ with $i_1 \leq j_1 \leq i_1 - r$ and $i_2 \leq j_2 \leq i_2 - r$.
\end{remark}

The 2-D system representation of convolutional Neural Networks facilitates the application of control theory methods to the analysis of CNNs. This includes finding an upper bound on the Lipschitz constant of CNNs. To this end, it is helpful that the Lipschitz constant falls together with the incremental $\ell_2$-gain of a system, which plays an important role in stability analysis.

\begin{definition}[The $\ell_2$-gain of 2-D systems]
    The $\ell_2$-gain of a 2-D system $S:\ell_2^{c_1}(\bbN\times \bbN) \to \ell_2^{c_2}(\bbN\times \bbN), u\lsb i_1,i_2 \rsb  \mapsto y\lsb i_1,i_2 \rsb $ is defined as the infimum over all constants $\gamma$, such that there exists a $c \in \bbR$ with
    \begin{align*}
        \left\| (y\lsb i_1,i_2 \rsb )\right\|_{\ell_2} &\leq c + \gamma \left\| (u\lsb i_1,i_2 \rsb )\right\|_{\ell_2} & \forall (u\lsb i_1,i_2 \rsb ) \in \ell_2, (y\lsb i_1,i_2 \rsb ) = S((u\lsb i_1,i_2 \rsb )).
    \end{align*}
    The same definition applies for the $\ell_2$-gain of a 2-D Lur´e system.
\end{definition}

To use the $\ell_2$-gain for finding a bound on the Lipschitz constant, we need to study the mapping from the difference of two input signals $(u^2\lsb i_1,i_2 \rsb  - u^1\lsb i_1,i_2 \rsb )$ to the respective difference  of output signals $(y^2\lsb i_1,i_2 \rsb  - y^1\lsb i_1,i_2 \rsb )$. For this mapping, the Lipschitz constant falls together with the $\ell_2$-gain. The following lemma shows that this mapping can be realized as a non-affine 2-D Lur´e system.

\begin{lemma}[Error dynamics of affine 2-D Lur´e systems]
    \label{lem:errorDynamics}
    Let an affine 2-D Lur´e system, defined by $f_1$, $f_2$, $g_1$, $g_2$, $\calA_{11}$, $\calA_{12}$, $\calA_{21}$, $\calA_{22}$, $\calB_{11}$, $\calB_{12}$, $\calB_{21}$,  $\calB_{22}$, $\calC_{11}$, $\calC_{12}$, $\calC_{21}$, $\calC_{22}$, $\calD_{11}$, $\calD_{12}$, $\calD_{21}$ and $\calD_{22}$, be given and consider two input signals $(u^1\lsb i_1,i_2 \rsb )$ and $(u^2\lsb i_1,i_2 \rsb )$ and corresponding output signals $(y^1\lsb i_1,i_2 \rsb )$ and $(y^2\lsb i_1,i_2 \rsb )$.
	Then there exists a non-affine 2-D Lur´e system mapping the input signal $(\tilde{u}\lsb i_1,i_2 \rsb) = (u^2\lsb i_1,i_2 \rsb  - u^1\lsb i_1,i_2 \rsb )$ to the output signal $(\tilde{y}\lsb i_1,i_2 \rsb ) = (y^2\lsb i_1,i_2 \rsb  - y^1\lsb i_1,i_2 \rsb )$. Such a Lur´e system is given by
	\begin{align}
		 \begin{pmatrix}
			e_1\lsb i_1+1,i_2 \rsb \\
			e_2\lsb i_1,i_2+1 \rsb \\
			\tilde{z}\lsb i_1,i_2 \rsb \\
			\tilde{y}\lsb i_1,i_2 \rsb 
		\end{pmatrix}
		=
		\begin{pmatrix}
			\calA_{11} & \calA_{12} & \calB_{11} & \calB_{12} \\
			\calA_{21} & \calA_{22} & \calB_{12} & \calB_{22} \\
			\calC_{11} & \calC_{12} & \calD_{11} & \calD_{12}\\
			\calC_{21} & \calC_{22} & \calD_{21} & \calD_{22}
		\end{pmatrix}
		\begin{pmatrix}
			e_1\lsb i_1-1,i_2 \rsb \\
			e_2\lsb i_1,i_2-1 \rsb \\
			\tilde{w}\lsb i_1,i_2 \rsb \\
			\tilde{u}\lsb i_1,i_2 \rsb 
		\end{pmatrix}, \hspace{5mm} \tilde{w}\lsb i_1,i_2 \rsb  = \tilde{\phi}_{i_1,i_2} (\tilde{z}\lsb i_1,i_2 \rsb ), \label{eq:error_2Dsys}
	\end{align}
	where the family of nonlinearities $\tilde{\phi}$ parametrized by $z^1 \in \ell_{2e}(\bbN \times \bbN)$ is defined by
	\begin{align*}
		\tilde{\phi}_{i_1,i_2}(\tilde{z}) &:= \phi (\tilde{z}  + z^1\lsb i_1,i_2 \rsb ) - \phi (z^1\lsb i_1,i_2 \rsb )
	\end{align*}
	for a real vector $\tilde{z}$ and $z^1\lsb i_1,i_2 \rsb $ is the input to $\phi$ of the original 2-D system for the input $u^1\lsb i_1,i_2 \rsb $.
\end{lemma}
\begin{proof}
	Consider the signals
	\begin{align*}
		e_1\lsb i_1,i_2 \rsb  &:= x^2_1\lsb i_1,i_2 \rsb  - x^1_1\lsb i_1,i_2 \rsb , \quad e_2\lsb i_1,i_2 \rsb  := x^2_2\lsb i_1,i_2 \rsb  - x^1_2\lsb i_1,i_2 \rsb , \quad \tilde{w}\lsb i_1,i_2 \rsb  := w^2\lsb i_1,i_2 \rsb  - w^1\lsb i_1,i_2 \rsb ,\\
		\tilde{z}\lsb i_1,i_2 \rsb  &:= z^2\lsb i_1,i_2 \rsb  - z^1\lsb i_1,i_2 \rsb , \quad \tilde{u}\lsb i_1,i_2 \rsb  := u^2\lsb i_1,i_2 \rsb  - u^1\lsb i_1,i_2 \rsb ,\quad \tilde{y}\lsb i_1,i_2 \rsb  := y^2\lsb i_1,i_2 \rsb  - y^1\lsb i_1,i_2 \rsb .
	\end{align*}
	We show that these signals satisfy the equations in \eqref{eq:error_2Dsys}. Then the system \eqref{eq:error_2Dsys} maps $\tilde{u}\lsb i_1,i_2 \rsb  = u^2\lsb i_1,i_2 \rsb  - u^1\lsb i_1,i_2 \rsb $ to $\tilde{y}\lsb i_1,i_2 \rsb  = y^2\lsb i_1,i_2 \rsb  - y^1\lsb i_1,i_2 \rsb $ as we claim.
	
	By the definition of the signals, we have that
	\begin{align*}
		\begin{pmatrix}
			e_1\lsb i_1,i_2 \rsb \\
			e_2\lsb i_1,i_2 \rsb \\
			\tilde{z}\lsb i_1,i_2 \rsb \\
			\tilde{y}\lsb i_1,i_2 \rsb 
		\end{pmatrix}
		&=
		\begin{pmatrix}
			x^2_1\lsb i_1,i_2 \rsb \\
			x^2_2\lsb i_1,i_2 \rsb \\
			z^2\lsb i_1,i_2 \rsb \\
			y^2\lsb i_1,i_2 \rsb 
		\end{pmatrix}
		-
		\begin{pmatrix}
			x^1_1\lsb i_1,i_2 \rsb \\
			x^1_2\lsb i_1,i_2 \rsb \\
			z^1\lsb i_1,i_2 \rsb \\
			y^1\lsb i_1,i_2 \rsb 
		\end{pmatrix}.
	\end{align*}
	Since $x^2_1, x_2^2,z^2$ and $y^2$ and $x_1^1, x_2^1,z^1$ and $y^1$ are solutions of System \eqref{eq:2D_linearFractionalRepresentation}, we get
	\begin{align*}
		\begin{pmatrix}
			e_1\lsb i_1+1,i_2 \rsb \\
			e_2\lsb i_1,i_2+1 \rsb \\
			\tilde{z}\lsb i_1,i_2 \rsb \\
			\tilde{y}\lsb i_1,i_2 \rsb 
		\end{pmatrix}
		&=
		\begin{pmatrix}
			f_1\\
			f_2\\
			g_1\\
			g_2
		\end{pmatrix}
		+
		\begin{pmatrix}
			\calA_{11} & \calA_{12} & \calB_{11} & \calB_{12} \\
			\calA_{21} & \calA_{22} & \calB_{12} & \calB_{22} \\
			\calC_{11} & \calC_{12} & \calD_{11} & \calD_{12}\\
			\calC_{21} & \calC_{22} & \calD_{21} & \calD_{22}
		\end{pmatrix}
		\begin{pmatrix}
			x_1^2\lsb i_1,i_2 \rsb \\
			x_2^2\lsb i_1,i_2 \rsb \\
			w^2\lsb i_1,i_2 \rsb \\
			u^2\lsb i_1,i_2 \rsb 
		\end{pmatrix}
		-
		\begin{pmatrix}
			f_1\\
			f_2\\
			g_1\\
			g_2
		\end{pmatrix}
		-
		\begin{pmatrix}
			\calA_{11} & \calA_{12} & \calB_{11} & \calB_{12} \\
			\calA_{21} & \calA_{22} & \calB_{12} & \calB_{22} \\
			\calC_{11} & \calC_{12} & \calD_{11} & \calD_{12}\\
			\calC_{21} & \calC_{22} & \calD_{21} & \calD_{22}
		\end{pmatrix}
		\begin{pmatrix}
			x_1^1\lsb i_1,i_2 \rsb \\
			x_2^1\lsb i_1,i_2 \rsb \\
			w^1\lsb i_1,i_2 \rsb \\
			u^1\lsb i_1,i_2 \rsb 
		\end{pmatrix}\\
		&=
		\begin{pmatrix}
			\calA_{11} & \calA_{12} & \calB_{11} & \calB_{12} \\
			\calA_{21} & \calA_{22} & \calB_{12} & \calB_{22} \\
			\calC_{11} & \calC_{12} & \calD_{11} & \calD_{12}\\
			\calC_{21} & \calC_{22} & \calD_{21} & \calD_{22}
		\end{pmatrix}
		\begin{pmatrix}
			x_1^2\lsb i_1,i_2 \rsb  - x_1^1\lsb i_1,i_2 \rsb \\
			x_2^2\lsb i_1,i_2 \rsb  - x_2^1\lsb i_1,i_2 \rsb \\
			w^2\lsb i_1,i_2 \rsb  - w^1\lsb i_1,i_2 \rsb \\
			u^2\lsb i_1,i_2 \rsb  - u^1\lsb i_1,i_2 \rsb 
		\end{pmatrix}\\
		&=
		\begin{pmatrix}
			\calA_{11} & \calA_{12} & \calB_{11} & \calB_{12} \\
			\calA_{21} & \calA_{22} & \calB_{12} & \calB_{22} \\
			\calC_{11} & \calC_{12} & \calD_{11} & \calD_{12}\\
			\calC_{21} & \calC_{22} & \calD_{21} & \calD_{22}
		\end{pmatrix}
		\begin{pmatrix}
			e_1\lsb i_1,i_2 \rsb \\
			e_2\lsb i_1,i_2 \rsb \\
			\tilde{w}\lsb i_1,i_2 \rsb \\
			\tilde{u}\lsb i_1,i_2 \rsb 
		\end{pmatrix}.
	\end{align*}
	It is left to show that the equation for the nonlinearity is satisfied. To this end, we consider
	\begin{align*}
		\tilde{w}\lsb i_1,i_2 \rsb  &= w^2\lsb i_1,i_2 \rsb  - w^1\lsb i_1,i_2 \rsb 
	\end{align*}
	Now, since $w^2\lsb i_1,i_2 \rsb  = \phi (z^2 \lsb i_1,i_2 \rsb )$ and $w^1\lsb i_1,i_2 \rsb  = \phi (z^1 \lsb i_1,i_2 \rsb )$ we obtain
	\begin{align*}
		\tilde{w}\lsb i_1,i_2 \rsb  &= \phi (z^2 \lsb i_1,i_2 \rsb ) - \phi (z^1 \lsb i_1,i_2 \rsb )\\
		&= \phi (\tilde{z} \lsb i_1,i_2 \rsb  - z^1 \lsb i_1,i_2 \rsb ) - \phi (z^1 \lsb i_1,i_2 \rsb )
		= \tilde{\phi}_{i_1,i_2} (\tilde{z} \lsb i_1,i_2 \rsb ),
	\end{align*}
	where we used $\tilde{z}\lsb i_1,i_2 \rsb  = z^2\lsb i_1,i_2 \rsb  - z^1\lsb i_1,i_2 \rsb $. This concludes the proof.
\end{proof}

Clearly, this representation has the disadvantage that the nonlinear function $\tilde{\phi}_{i_1,i_2}$ depends on the signal $z^1\lsb i_1,i_2 \rsb $. For this reason, it is common in the literature to derive a joint bound on the $\ell_2$-gain of \eqref{eq:error_2Dsys}, that applies for all possible versions of $\tilde{\phi}_{i_1,i_2}(\cdot)$, i.e., for all signals $(z^1 \lsb i_1,i_2 \rsb)$ parametrizing this nonlinearity.

\begin{corollary}
	\label{cor:errorDynamicGainBound}
	Assume that the 2-D Lur´e system \eqref{eq:error_2Dsys} satisfies the $\ell_2$-gain bound
	\begin{align*}
		\|(\tilde{y}\lsb i_1,i_2 \rsb )\|_{\ell_2} \leq \gamma \|(\tilde{u}\lsb i_1,i_2 \rsb )\|_{\ell_2}
	\end{align*}
	for all nonlinearities of the form $\tilde{z} \mapsto \tilde{\phi}_{i_1,i_2}(\tilde{z}) := \phi (\tilde{z}\lsb i_1,i_2 \rsb  + z^1\lsb i_1,i_2 \rsb ) - \phi (z^1\lsb i_1,i_2 \rsb )$, where $z^1\lsb i_1,i_2 \rsb $ is any signal in $\ell_{2e}(\bbN \times \bbN)$. Then the 2-D system \eqref{eq:2D_linearFractionalRepresentation} has an incremental $\ell_2$-gain (Lipschitz constant) of at most $\gamma$.
\end{corollary}

\begin{proof}
	Trivial.
\end{proof}

Due to the result of Corollary \ref{cor:errorDynamicGainBound}, we are left with the problem of finding the $\ell_2$-gain of a 2-D Lur´e system for all nonlinearities of the form $\tilde{\phi}_{i_1,i_2}(\tilde{z}\lsb i_1,i_2 \rsb ) := \phi (\tilde{z}\lsb i_1,i_2 \rsb  + z^1\lsb i_1,i_2 \rsb ) - \phi (z^1\lsb i_1,i_2 \rsb )$, where $z^1\lsb i_1,i_2 \rsb $ is any signal in $\ell_{2e}(\bbN\times \bbN)$. This problem is addressed in the next section by means of a dissipativity theory for 2-D systems.

\section{Dissipativity and IQCs for 2-D systems}
\label{sec:3}

The layers of convolutional Neural Networks can be realized as Roesser systems (Definition \ref{def:RoesserSystem}). These Roesser systems constitute linear shift-invariant functions $\ell_{2e}^{c_1}(\bbN \times \bbN) \to \ell_{2e}^{c_2}(\bbN \times \bbN), (u\lsb i_1,i_2 \rsb ) \mapsto (y\lsb i_1,i_2 \rsb )$ with the state space realization
\begin{align}
    \begin{split}
    	x_1\lsb i_1+1,i_2 \rsb  &= \bA_{11} x_1\lsb i_1,i_2 \rsb  + \bA_{12} x_2\lsb i_1,i_2 \rsb  + \bB_1 u\lsb i_1,i_2 \rsb\\
    x_2\lsb i_1,i_2+1 \rsb  &= \bA_{21} x_1\lsb i_1,i_2 \rsb  + \bA_{22} x_2\lsb i_1,i_2 \rsb  + \bB_2 u\lsb i_1,i_2 \rsb\\
    y\lsb i_1,i_2 \rsb  &= \bC_1 x_1\lsb i_1,i_2 \rsb  + \bC_2 x_2\lsb i_1,i_2 \rsb  + \bD u\lsb i_1,i_2 \rsb.
    \end{split}
    \label{eq:2DsysDissi}
\end{align}
The key difference to 1-D systems is that states, inputs, and outputs depend on two indices $i_1$ and $i_2$ and that information is transmitted along two directions - from $i_1$ to $i_1 + 1$ and from $i_2$ to $i_2 + 1$.

A key role in the system analysis for convolutional Neural Networks and 2-D systems proposed in this work is played by dissipativity theory.
This control theory tool, originally proposed for 1-D systems \cite{willems1972dissipative1,willems1972dissipative2} and generalized to Roesser systems in \cite{ahn2015two} extends Lyapunov theory to systems with inputs and outputs and is a cornerstone of robust control theory. Indeed, the notion of supply rates from dissipativity theory enables the specification of important control objectives (such as a bound on the $\ell_2$ gain) and robustness properties. Celebrated tools from robust control theory, such as integral quadratic constraints can be integrated into the dissipativity framework. Also, the original works on estimating the Lipschitz constant of fully connected Neural Networks, \cite{fazlyab2019efficient}, can be viewed as an application of dissipativity theory.
To transfer dissipativity theory to 2-D systems, we use the following notion of dissipativity (compare \cite{ahn2015two}).

\begin{definition}[Dissipativity]
    \label{def:dissipativity}
    We call a 2-D system \eqref{eq:2DsysDissi} dissipative with respect to a supply rate $s: \bbR^{c_1} \times \bbR^{c_2} \to \bbR, (u,y)\mapsto s(u,y)$, if there exist two storage functions $V_1:\bbR^{n_1} \to \bbR_{\geq 0}$ and $V_2:\bbR^{n_2}\to \bbR_{\geq 0}$ such that any input output state trajectory $(x_1\lsb i_1,i_2 \rsb ,x_2\lsb i_1,i_2 \rsb ,u\lsb i_1,i_2 \rsb ,y\lsb i_1,i_2 \rsb )$ of \eqref{eq:2DsysDissi} satisfies for all $N_1,N_2 \in \bbN$
    \begin{align}
        \sum_{i_2 = 0}^{N_2-1} \left(V_1(x_1\lsb N_1,i_2 \rsb ) - V_1(x_1\lsb 0,i_2 \rsb )\right)+ \sum_{i_1 = 0}^{N_1-1} \left(V_2(x_2\lsb i_1,N_2 \rsb ) - V_2(x_2\lsb i_1,0 \rsb )\right) \leq \sum_{i_1 = 0}^{N_1-1} \sum_{i_2 = 0}^{N_2-1} s(u\lsb i_1,i_2 \rsb,y\lsb i_1,i_2 \rsb ). \label{eq:dissipativity_condition}
    \end{align}
\end{definition}


Note that in analogy to classical dissipativity for 1-D systems, there is an energy interpretation to the storage functions $V_1$ and $V_2$ and to the supply $s$. The functions $V_1$ and $V_2$ measure the energy that is transported from node $\lsb i_1,i_2-1 \rsb $ and from node $\lsb i_1-1,i_2 \rsb $ to node $\lsb i_1,i_2 \rsb $ and $s$ measures how much energy is supplied to the system at node $\lsb i_1,i_2 \rsb $. The sum of these energies must be larger than or equal to the sum of the energies that node $\lsb i_1,i_2 \rsb $ transfers to $\lsb i_1+1,i_2 \rsb $ and $\lsb i_1,i_2+1 \rsb $.  Stated otherwise, the \emph{energy divergence} at node $\lsb i_1,i_2 \rsb $ must be bounded by the energy supply at node $\lsb i_1,i_2 \rsb $.

This definition of dissipativity generalizes the common dissipativity notion for 1-D systems. For the usefulness of dissipativity theory, it is important that we can easily verify whether a system is dissipative for given storage functions and supply rates. This is ensured by dissipation inequalities.

\begin{theorem}[Dissipation inequality]
	\label{thm:dissipativity}
    A 2-D system \eqref{eq:2DsysDissi} is dissipative according to Definition \ref{def:dissipativity} if and only if the dissipation inequality
    \begin{align}\label{eq:dis_inequality}
        V_1(x_1\lsb i_1+1,i_2 \rsb ) + V_2(x_2\lsb i_1,i_2+1 \rsb ) &\leq V_1(x_1\lsb i_1,i_2 \rsb ) + V_2(x_2\lsb i_1,i_2 \rsb ) + s(u\lsb i_1,i_2 \rsb,y\lsb i_1,i_2 \rsb )
    \end{align}
    is satisfied for all $i_1,i_2 \in \bbN$ and for all admissible system trajectories $x_1,x_2,u,y$.
    The dissipation inequality implies
    \begin{align}
        \sum_{i_1 = 0}^{N_1-1} \sum_{i_2 = 0}^{N_2-1} s(u\lsb i_1,i_2 \rsb,y\lsb i_1,i_2 \rsb ) &\geq 0 & \forall N_1,N_2 \in \bbN \label{eq:dis-inequality2}
    \end{align}
    if $V_1(x_1\lsb i_1,0 \rsb ) = 0~ \forall i_1 \in \bbN$ and if $V_2(x_2\lsb 0,i_2 \rsb ) = 0 ~ \forall i_2 \in \bbN$. 
\end{theorem}

\begin{proof}
    The dissipation inequality is a special case of \eqref{eq:dissipativity_condition} for $N_1 = N_2 = 1$. Consequently, all that is to show is that the dissipation inequality \eqref{eq:dis_inequality} implies \eqref{eq:dissipativity_condition}. To this end, we sum \eqref{eq:dis_inequality} for all $i_1 = 0,\ldots , N_1 - 1$ and $i_2 = 0 ,\ldots, N_2 - 1$ yielding
    \begin{align*}
        \sum_{i_1 = 0}^{N_1 - 1} \sum_{i_2 = 0}^{N_2 - 1} s(u\lsb i_1,i_2 \rsb,y\lsb i_1,i_2 \rsb ) &\geq 
        \sum_{i_2 = 0}^{N_2 - 1} \sum_{i_1 = 0}^{N_1 - 1} \left(V_1(x_1\lsb i_1+1,i_2 \rsb ) - V_1(x_1\lsb i_1,i_2 \rsb )\right)\\
        &+ \sum_{i_1 = 0}^{N_1 - 1} \sum_{i_2 = 0}^{N_2 - 1} \left( V_2(x_2\lsb i_1,i_2+1 \rsb ) - V_2(x_2\lsb i_1,i_2 \rsb )\right).
    \end{align*}
    By replacing the telescopic sums $\sum_{i_1 = 0}^{N_1 - 1} \left(V_1(x_1\lsb i_1+1,i_2 \rsb ) - V_1(x_1\lsb i_1,i_2 \rsb )\right)$ and $\sum_{i_2 = 0}^{N_2 - 1} \left( V_2(x_2\lsb i_1,i_2+1 \rsb ) - V_2(x_2\lsb i_1,i_2 \rsb )\right)$ with $V_1(x_1\lsb N_1,i_2 \rsb ) - V_1(x_1\lsb 0,i_2 \rsb )$ and $V_2(x_2\lsb i_1,N_2 \rsb ) - V_2(x_2\lsb i_1,0 \rsb )$ respectively, we obtain \eqref{eq:dissipativity_condition}.
    That \eqref{eq:dis-inequality2} is implied by \eqref{eq:dis_inequality} can then simply be seen by setting $V_1(x_1\lsb 0,i_1 \rsb )$ and $V_2(x_2\lsb i_2,0 \rsb )$ equal to zero in \eqref{eq:dissipativity_condition} and observing that $V_1(x_1\lsb N_1,i_2 \rsb )$ and $V_2(x_2\lsb i_1,N_2 \rsb )$ are always non-negative.
\end{proof}

It should be stressed that dissipativity with respect to supply rates $s(u,y)$ includes a certificate for an $\ell_2$-gain bound of $\gamma$ as a special case. Indeed, by choosing the supply rate
\begin{align}
    s(u,y) := \begin{pmatrix}
        u\\ y
    \end{pmatrix}^\top
    \begin{pmatrix}
        \gamma^2 \bI & 0\\
        0 & -\bI
    \end{pmatrix}
    \begin{pmatrix}
        u\\ y
    \end{pmatrix}, \label{eq:LipschitzSupply}
\end{align}
we obtain for $(u\lsb i_1,i_2 \rsb ) \in \ell_2(\bbN \times \bbN)$ and $V_1(x_1\lsb i_1,0 \rsb ) = 0 ~ \forall i_1 \in \bbN$ and $V_2(x_2\lsb 0,i_2 \rsb ) = 0 ~ \forall i_2 \in \bbN$ from Theorem \ref{thm:dissipativity} that
\begin{align*}
    0 \leq \sum_{i_1,i_2 = 0}^\infty s(u\lsb i_1,i_2 \rsb ,y\lsb i_1,i_2 \rsb ) = \sum_{i_1,i_2 = 0}^\infty (\|y\lsb i_1,i_2 \rsb \|^2_2 - \gamma^2 \|u\lsb i_1,i_2 \rsb \|_2^2) = \gamma^2 \|(u\lsb i_1,i_2 \rsb )\|_{\ell_2}^2 - \|(y\lsb i_1,i_2 \rsb )\|_{\ell_2}^2
\end{align*}
holds, which is equivalent to an $\ell_2$-gain of $\gamma$ for the considered system.

Next, we consider robust dissipativity. Apart from being a useful generalization of dissipativity in general, this notion is going to allow us to deal with the $\ell_2$-gain problem from Corollary \ref{cor:errorDynamicGainBound}, where dissipativity needs to be satisfied for all $\tilde{\phi}_{i_1,i_2}(\tilde{z}\lsb i_1,i_2 \rsb ) := \phi (\tilde{z}\lsb i_1,i_2 \rsb  + z^1\lsb i_1,i_2 \rsb ) - \phi (z^1\lsb i_1,i_2 \rsb )$, where $z^1\lsb i_1,i_2 \rsb $ is any signal in $\ell_{2e}(\bbN\times \bbN)$.

\begin{theorem}[Robust dissipativity for Lur´e systems]
    \label{thm:robustDissipativity}
    Consider a well-posed 2-D Lur´e system \eqref{eq:2D_linearFractionalRepresentation}
    where the nonlinear component $\phi$ is known to satisfy the inequality
    \begin{align}
        s^w(z,\phi(z)) & \leq 0 \qquad \forall z \in \bbR^{m}. \label{eq:dissiUncertainty}
    \end{align}
    Then this Lur´e system is dissipative with respect to the supply rate $s: (u,y) \mapsto s(u,y)$, i.e., all input-output pairs $u,y$ of this Lur´e system satisfy \eqref{eq:dissipativity_condition}, if there exist storage functions $V_1: \bbR^{n_1} \to \bbR_{\geq 0 }$ and $V_2: \bbR^{n_2} \to \bbR_{\geq 0}$ satisfying the robust dissipation inequality
    \begin{align}
        V_1(x_1\lsb i_1+1,i_2 \rsb ) + V_2(x_2\lsb i_1,i_2+1 \rsb ) 
        &\leq V_1(x_1\lsb i_1,i_2 \rsb ) + V_2(x_2\lsb i_1,i_2 \rsb ) + s(u\lsb i_1,i_2 \rsb,y\lsb i_1,i_2 \rsb ) + s^w(z\lsb i_1,i_2 \rsb ,w\lsb i_1,i_2 \rsb ) \label{eq:dissipativityRobust}
    \end{align}
	for all $i_1,i_2 \in \bbN$ and all trajectories $x_1\lsb i_1,i_2 \rsb $, $x_2\lsb i_1,i_2 \rsb $, $u\lsb i_1,i_2 \rsb $, $w\lsb i_1,i_2 \rsb $, $z\lsb i_1,i_2 \rsb $, $y\lsb i_1,i_2 \rsb $ of \eqref{eq:2D_linearFractionalRepresentation}.
\end{theorem}

\begin{proof}
	To prove that the 2-D Lur´e system is dissipative, it suffices to show that for all admissible trajectories consisting of $x_1\lsb i_1,i_2 \rsb $, $x_2\lsb i_1,i_2 \rsb $, $u\lsb i_1,i_2 \rsb $, $w\lsb i_1,i_2 \rsb $, $z\lsb i_1,i_2 \rsb $, $y\lsb i_1,i_2 \rsb $, the dissipation inequality
	\begin{align*}
		V_1(x_1\lsb i_1+1,i_2 \rsb ) + V_2(x_2\lsb i_1,i_2+1 \rsb ) &\leq V_1(x_1\lsb i_1,i_2 \rsb ) + V_2(x_2\lsb i_1,i_2 \rsb ) + s(u\lsb i_1,i_2 \rsb,y\lsb i_1,i_2 \rsb )
	\end{align*}
	is satisfied. A fact, which follows directly from \eqref{eq:dissipativityRobust} by taking into account that the trajectory $x_1,x_2,u,y,z,w$ satisfies $w\lsb i_1,i_2 \rsb  = \phi (w\lsb i_1,i_2 \rsb )$ for all $i_1,i_2 \in \bbN$, which in turn implies $s^w(z\lsb i_1,i_2 \rsb ,w\lsb i_1,i_2 \rsb ) \leq 0$.
\end{proof}

The proof of Theorem \ref{thm:robustDissipativity} also allows for an energy interpretation. Namely, since $s^w(z,\phi(z))$ is always smaller than or equal to zero, the nonlinearity can only take energy away from the system. Consequently, it can only make the dissipativity with respect to the supply rate $(u,y) \mapsto s(u,y)$ more strict.
The key advantage of \eqref{eq:dissipativityRobust} is that the nonlinearity $\phi$ does not appear in the dissipation inequality anymore which often makes it easier to verify \eqref{eq:dissipativityRobust} instead of \eqref{eq:dissipativity_condition} with a given $\phi(\cdot)$ (or possibly a whole family of nonlinearities, which we need to consider for Neural Networks).

Besides applications in robust control, it turns out that robust dissipativity is a useful tool to estimate the Lipschitz constant of Neural Networks (see e.g. \cite{fazlyab2019efficient}). This is because the typical activation functions $\sigma(\cdot)$ of Neural Networks, such as ReLU, Sigmoid, tanh, ELU, and many others all satisfy
\begin{align*}
    \frac{\sigma (t_1) - \sigma (t_2)}{t_1 - t_2} \in [0,1] \quad \forall t_1\neq t_2 \in \bbR,
\end{align*}
i.e., their slope is restricted to the interval $[0,1]$. Slope-restricted nonlinearities are long studied in robust control. Indeed, the slope restriction implies that the uncertainty $\phi = \sigma$ satisfies the condition
\begin{equation}
        0\geq s^w(z_2 - z_1,\sigma (z_2) - \sigma (z_1))
    =
    \begin{pmatrix}
        \sigma (z_2) - \sigma (z_1)\\
        z_2 - z_1
    \end{pmatrix}^\top
    \begin{pmatrix}
        2\Lambda & -\Lambda\\
        -\Lambda & 0
    \end{pmatrix}
    \begin{pmatrix}
        \sigma (z_2) - \sigma (z_1)\\
        z_2 - z_1
    \end{pmatrix}
    \qquad \forall z \in \bbR^m
    \label{eq:sectorZeroOneSupply}
\end{equation}
for all matrices $\Lambda = \diag \lambda$ with $\lambda \in \bbR^{m}_{\geq 0}$ and $m \in \bbN$, i.e., diagonal matrices with non-negative entries. Note that if $\phi = \sigma$ in Lemma \ref{lem:errorDynamics}, then we obtain $\tilde{\phi}_{i_1,i_2}(\tilde{z}\lsb i_1,i_2 \rsb ) = \sigma (z^2\lsb i_1,i_2 \rsb) - \sigma (z^1\lsb i_1,i_2 \rsb)$ and $\tilde{z}\lsb i_1,i_2 \rsb = z^2\lsb i_1,i_2 \rsb - z^1\lsb i_1,i_2 \rsb$, implying that $s^w(\tilde{z}\lsb i_1,i_2 \rsb ,\tilde{\phi}_{i_1,i_2}(\tilde{z}\lsb i_1,i_2 \rsb )) \leq 0$ holds also true. Combined with the fact that our objective to find an upper bound on the $\ell_2$-gain can be formulated with the quadratic supply rate \eqref{eq:sectorZeroOneSupply}, we can formulate the conditions of Corollary \ref{cor:errorDynamicGainBound} as a robust dissipativity problem. For this reason, we will be dealing with dissipativity of 2-D systems with respect to quadratic supply rates in the remainder of this section. In addition to the supply rate, we also choose the storage function quadratic as it is common practice for 1-D systems, where it is also shown that quadratic storage functions are not a limiting assumption for dissipativity with respect to quadratic supply rates. Hence, we consider the storage functions and supply rates
\begin{align}
    s(u,y) &= \begin{pmatrix}
        u\\ y
    \end{pmatrix}^\top
    \begin{pmatrix}
        \bR & \bS\\
        \bS^\top & \bQ
    \end{pmatrix}
    \begin{pmatrix}
        u\\ y
    \end{pmatrix}, &
    V_1(x_1) &= x^\top_1 \bP_1 x_1, & V_2(x_2) &= x^\top_2 \bP_2 x_2. \label{eq:quadraticSupplyStorage}
\end{align}
To verify the dissipativity of a 2-D system with respect to the supply rate and storage functions as above, we utilize the dissipation inequality, which leads us to the problem of verifying the following property. We need to find $\bP_1\succeq 0$ and $\bP_2 \succeq 0$, such that
\begin{align*}
    x_1\lsb i_1+1,i_2 \rsb ^\top \bP_1 x_1\lsb i_1+1,i_2 \rsb  + x_2\lsb i_1,i_2+1 \rsb ^\top \bP_2 x_2\lsb i_1,i_2+1 \rsb  &\leq x_1\lsb i_1,i_2 \rsb  \bP_1 x_1\lsb i_1,i_2 \rsb  + x_2\lsb i_1,i_2 \rsb ^\top \bP_2 x_2\lsb i_1,i_2 \rsb \\
    &+
    \begin{pmatrix}
        u\lsb i_1,i_2 \rsb\\ y\lsb i_1,i_2 \rsb 
    \end{pmatrix}^\top
    \begin{pmatrix}
        \bR & \bS\\
        \bS^\top & \bQ
    \end{pmatrix}
    \begin{pmatrix}
        u\lsb i_1,i_2 \rsb\\ y\lsb i_1,i_2 \rsb 
    \end{pmatrix}
\end{align*}
is satisfied for all possible values of $x_1\lsb i_1,i_2 \rsb ,x_1\lsb i_1+1,i_2 \rsb ,x_2\lsb i_1,i_2 \rsb ,x_2\lsb i_1,i_2+1 \rsb ,u\lsb i_1,i_2 \rsb ,y\lsb i_1,i_2 \rsb $ satisfying the dynamic equations of \eqref{eq:2DsysDissi}.
At first glance, this might seem like a hard problem, since an infinite number of values for the signals needs to be checked for this condition. Luckily, in robust control, powerful linear matrix inequality (LMI) techniques have emerged to deal with such dissipativity problems. The following theorem shows, that also for 2-D systems, dissipativity can be certified using linear matrix inequalities.

\begin{theorem}[Quadratic dissipativity certificate]
    \label{thm:5}
    The 2-D system \eqref{eq:2DsysDissi} is dissipative with respect to the supply rate \eqref{eq:quadraticSupplyStorage} if there exist $\bP_1 \succeq 0$ and $\bP_2 \succeq 0$ satisfying the linear matrix inequality
    \begin{align}
        \begin{pmatrix}
        	\bA_{11} & \bA_{12} & \bB_1\\
            \bA_{21} & \bA_{22} & \bB_2\\
            \bI & 0 & 0\\
            0 & \bI & 0\\
            \bC_1 & \bC_2 & \bD\\
            0 & 0 & \bI
        \end{pmatrix}^\top
        \left(
        \begin{array}{cccc|cc}
        	-\bP_1 & & & &\\
            & -\bP_2 & & &\\
            & & \bP_1 & & &\\
            & & & \bP_2 & &\\\hline
            & & & & \bQ & \bS\\
            & & & & \bS^\top & \bR
        \end{array}
        \right)
        \begin{pmatrix}
        	\bA_{11} & \bA_{12} & \bB_1\\
        	\bA_{21} & \bA_{22} & \bB_2\\
        	\bI & 0 & 0\\
        	0 & \bI & 0\\
        	\bC_1 & \bC_2 & \bD\\
        	0 & 0 & \bI
        \end{pmatrix}
        \succeq 0. \label{eq:2DsysLMI}
    \end{align}
\end{theorem}

\begin{proof}
    This theorem is equivalent to Theorem 1 from \cite{ahn2015two}.
    Let $(x_1\lsb i_1,i_2 \rsb ), (x_2\lsb i_1,i_2 \rsb ) \in \ell_2(\bbN\times \bbN), (u\lsb i_1,i_2 \rsb ) \in \ell_2(\bbN\times \bbN)$ and $(y\lsb i_1,i_2 \rsb ) \in \ell_2(\bbN\times \bbN)$ be signals satisfying the system equations \eqref{eq:2DsysDissi}. Multiplying \eqref{eq:2DsysLMI} from left and right with $\begin{pmatrix}
		x\lsb i_1,i_2 \rsb ^\top & x\lsb i_1,i_2 \rsb ^\top & u\lsb i_1,i_2 \rsb ^\top
	\end{pmatrix}^\top$, we obtain
	\begin{align*}
		0 &\leq 
		\begin{pmatrix}
			x_1\lsb i_1,i_2 \rsb  \\ x_2\lsb i_1,i_2 \rsb  \\ u\lsb i_1,i_2 \rsb 
		\end{pmatrix}^\top
		\begin{pmatrix}
			\bA_{11} & \bA_{12} & \bB_1\\
			\bA_{21} & \bA_{22} & \bB_2\\
			\bI & 0 & 0\\
			0 & \bI & 0\\
			\bC_1 & \bC_2 & \bD\\
			0 & 0 & \bI
		\end{pmatrix}^\top
		\left(
		\begin{array}{cccc|cc}
			-\bP_1 & & & &\\
			& -\bP_2 & & &\\
			& & \bP_1 & & &\\
			& & & \bP_2 & &\\\hline
			& & & & \bQ & \bS\\
			& & & & \bS^\top & \bR
		\end{array}
		\right)
		\begin{pmatrix}
			\bA_{11} & \bA_{12} & \bB_1\\
			\bA_{21} & \bA_{22} & \bB_2\\
			\bI & 0 & 0\\
			0 & \bI & 0\\
			\bC_1 & \bC_2 & \bD\\
			0 & 0 & \bI
		\end{pmatrix}
		\begin{pmatrix}
			x_1\lsb i_1,i_2 \rsb  \\ x_2\lsb i_1,i_2 \rsb  \\ u\lsb i_1,i_2 \rsb 
		\end{pmatrix}\\
		&=
		\begin{pmatrix}
			x_1\lsb i_1+1,i_2 \rsb \\x_2\lsb i_1,i_2+1 \rsb \\x_1\lsb i_1,i_2 \rsb \\x_2\lsb i_1,i_2 \rsb \\ y_\lsb i_1,i_2 \rsb \\ u\lsb i_1,i_2 \rsb 
		\end{pmatrix}^\top
			\left(
		\begin{array}{cccc|cc}
		-\bP_1 & & & &\\
		& -\bP_2 & & &\\
		& & \bP_1 & & &\\
		& & & \bP_2 & &\\\hline
		& & & & \bQ & \bS\\
		& & & & \bS^\top & \bR
		\end{array}
		\right)
		\begin{pmatrix}
			x_1\lsb i_1+1,i_2 \rsb \\x_2\lsb i_1,i_2+1 \rsb \\x_1\lsb i_1,i_2 \rsb \\x_2\lsb i_1,i_2 \rsb \\ y_\lsb i_1,i_2 \rsb \\ u\lsb i_1,i_2 \rsb 
		\end{pmatrix}\\
		&=
		- x_1\lsb i_1+1,i_2 \rsb ^\top \bP_1 x_1\lsb i_1+1,i_2 \rsb  - x\lsb i_1,i_2+1 \rsb ^\top \bP_2 x_2\lsb i_1,i_2+1 \rsb  + x_1\lsb i_1,i_2 \rsb ^\top \bP_1 x_1\lsb i_1,i_2 \rsb \\
		&+ x_2\lsb i_1,i_2 \rsb ^\top \bP_2 x_2\lsb i_1,i_2 \rsb 
		+
		\begin{pmatrix}
			u\lsb i_1,i_2 \rsb \\ y\lsb i_1,i_2 \rsb 
		\end{pmatrix}^\top
		\begin{pmatrix}
			\bR & \bS\\
			\bS^\top & \bQ
		\end{pmatrix}
		\begin{pmatrix}
			u\lsb i_1,i_2 \rsb \\ y\lsb i_1,i_2 \rsb 
		\end{pmatrix},
	\end{align*}
	which is the dissipation inequality for the supply rate in \eqref{eq:quadraticSupplyStorage} and the storage functions $V_2(x_2) = x_2^\top \bP_2 x_2$ and $V_1(x_1) = x_1^\top \bP_1 x_1$.
\end{proof}

Note that the matrix inequality conditions of Theorem \ref{thm:5} are linear (convex) in the decision variables $\bP_1$ and $\bP_2$, such that we can solve for those conditions using standard off-the-shelf LMI solvers. Inequality \eqref{eq:2DsysLMI} can be used, for example, to find the $\ell_2$-gain of a 2-D system (convolutional layer) by choosing $\bQ = \bI$, $\bS = 0$ and $\bR = -\gamma^2 \bI$ and minimizing $\gamma^2$ as the objective function subject to the LMI constraints of Theorem \ref{thm:5} (where $\gamma^2$ is another decision variable).

To apply this approach to fully convolutional Neural Networks, Theorem \ref{thm:5} needs to be generalized to the setting of Lur´e systems/robust dissipativity. Here, our approach is to treat the nonlinearity $\phi$ in a robust fashion, i.e., by means of a supply rate $s^w(\phi(z),z)$ which is smaller than or equal to zero for all $z$. The respective result is formulated in the following theorem.
\begin{theorem}
    \label{thm:6}
    Consider a 2-D Lur´e system of the form \eqref{eq:2D_linearFractionalRepresentation} with a nonlinearity $\phi$ satisfying $s^w(\phi(z),z) \leq 0$ for all $z$. Assume further that $s$ and $s^w$ are quadratic supply rates as in \eqref{eq:quadraticSupplyStorage} with the matrices $(\bQ,\bS,\bR)$ describing $s$ and $(\bQ^w,\bS^w,\bR^w)$ describing $s^w$. Then this system is dissipative with respect to the supply rate $(u,y) \mapsto s(u,y)$, if there exist $\bP_1 \succeq 0$ and $\bP_2 \succeq 0$, satisfying the matrix inequality
    \begin{align}
        \begin{pmatrix}
            \calA_{11} & \calA_{12} & \calB_{11} & \calB_{12}\\
            \calA_{12} & \calA_{22} & \calB_{21} & \calB_{22}\\
            \bI & 0 & 0 & 0\\
            0 & \bI & 0 & 0\\
            \calC_{11} & \calC_{12} & \calD_{11} & \calD_{12}\\
            0 & 0 & \bI & 0\\
            \calC_{21} & \calC_{22} & \calD_{21} & \calD_{22}\\
            0 & 0 & 0 & \bI
        \end{pmatrix}^\top
        \left(
        \begin{array}{cccc|cc|cc}
        	-\bP_1 & & & & &\\
            & -\bP_2 & & & &\\
            & &\bP_1 & & &\\
            & & & \bP_2 & &\\\hline
            & & & & \bQ & \bS\\
            & & & & \bS^\top & \bR\\\hline
            & & & & & & \bQ^w & \bS^w\\
            & & & & & & (\bS^w)^\top & \bR^w
        \end{array}
        \right)
        \begin{pmatrix}
        	\calA_{11} & \calA_{12} & \calB_{11} & \calB_{12}\\
        	\calA_{12} & \calA_{22} & \calB_{21} & \calB_{22}\\
        	\bI & 0 & 0 & 0\\
        	0 & \bI & 0 & 0\\
        	\calC_{11} & \calC_{12} & \calD_{11} & \calD_{12}\\
        	0 & 0 & \bI & 0\\
        	\calC_{21} & \calC_{22} & \calD_{21} & \calD_{22}\\
        	0 & 0 & 0 & \bI
        \end{pmatrix}
        \succeq 0. \label{eq:robustDissipativityLMI}
    \end{align}
\end{theorem}

\begin{proof}
	To prove this statement, we assume again, that a trajectory of the system \eqref{eq:2D_linearFractionalRepresentation} is given and multiply \eqref{eq:robustDissipativityLMI} from left and right by $\begin{pmatrix}
		x\lsb i_1,i_2 \rsb ^\top & x\lsb i_1,i_2 \rsb ^\top & w\lsb i_1,i_2 \rsb ^\top & u\lsb i_1,i_2 \rsb ^\top
	\end{pmatrix}^\top$ and its transpose. Due to
	\begin{align*}
		\begin{pmatrix}
			x_1\lsb i_1+1,i_2 \rsb \\
			x_2\lsb i_1,i_2+1 \rsb \\
			x_1\lsb i_1,i_2 \rsb \\
			x_2\lsb i_1,i_2 \rsb \\
			y\lsb i_1,i_2 \rsb \\
			u\lsb i_1,i_2 \rsb \\
			z\lsb i_1,i_2 \rsb \\
			w\lsb i_1,i_2 \rsb 
		\end{pmatrix}
		=
		\begin{pmatrix}
			\calA_{11} & \calA_{12} & \calB_{11} & \calB_{12}\\
			\calA_{12} & \calA_{22} & \calB_{21} & \calB_{22}\\
			\bI & 0 & 0 & 0\\
			0 & \bI & 0 & 0\\
			\calC_{11} & \calC_{12} & \calD_{11} & \calD_{12}\\
			0 & 0 & \bI & 0\\
			\calC_{21} & \calC_{22} & \calD_{21} & \calD_{22}\\
			0 & 0 & 0 & \bI
		\end{pmatrix}
		\begin{pmatrix}
			x_1\lsb i_1,i_2 \rsb \\
			x_2\lsb i_1,i_2 \rsb \\
			u\lsb i_1,i_2 \rsb \\
			w\lsb i_1,i_2 \rsb 
		\end{pmatrix},
	\end{align*}
	this yields
	\begin{align*}
		& -x_1\lsb i_1+1,i_2 \rsb ^\top \bP_1 x_1\lsb i_1+1,i_2 \rsb  - x\lsb i_1,i_2+1 \rsb ^\top \bP_2 x_2\lsb i_1,i_2+1 \rsb  + x_1\lsb i_1,i_2 \rsb ^\top \bP_1 x_1\lsb i_1,i_2 \rsb \\
		&+ x_2\lsb i_1,i_2 \rsb ^\top \bP_2 x_2\lsb i_1,i_2 \rsb 
		+
		\begin{pmatrix}
			u\lsb i_1,i_2 \rsb \\ y\lsb i_1,i_2 \rsb 
		\end{pmatrix}^\top
		\begin{pmatrix}
			\bR & \bS\\
			\bS^\top & \bQ
		\end{pmatrix}
		\begin{pmatrix}
			u\lsb i_1,i_2 \rsb \\ y\lsb i_1,i_2 \rsb 
		\end{pmatrix}
		+
		\begin{pmatrix}
			w\lsb i_1,i_2 \rsb \\ z\lsb i_1,i_2 \rsb 
		\end{pmatrix}^\top
		\begin{pmatrix}
			\bR^w & \bS^w\\
			(\bS^w)^\top & \bQ^w
		\end{pmatrix}
		\begin{pmatrix}
			w\lsb i_1,i_2 \rsb \\ z\lsb i_1,i_2 \rsb 
		\end{pmatrix} \leq 0,
	\end{align*}
	which is the robust dissipation inequality.
\end{proof}

Again, we emphasize that this matrix inequality condition is linear in $\bP_2$, $\bP_1$ and possible decision variables in $\bQ,\bS,\bR$ and $\bQ^w$, $\bS^w$, $\bR^w$.

\section{Convolutional Neural Networks with fully connected layers}
\label{sec:4}

In this section, we extend our fully convolutional Neural Network structure to include additional fully connected layers, commonly found in image classification. The architecture consists of a sequence of convolutional layers, followed by a sequence of fully connected layers alternately concatenated with activation function layers. A challenge for the analysis of these combined networks is the different input and output domains of the individual layers. Convolutional layers take 2-D signals as inputs and outputs, while fully connected layers take finite-dimensional vectors as inputs and outputs. To overcome this challenge, we introduce the embedding operator $\calE$ and the flattening operator $\calF$, and define our \emph{hybrid Neural Network} (HNN) as
\begin{align}
	\mathrm{HNN}_\theta = \calL_{l_2} \circ \sigma \circ \ldots \circ \sigma \circ \calL_{l+1} \circ \calF \circ \sigma \circ \calC_{l}\circ \sigma \circ \calC_{l - 1} \circ \sigma \circ \cdots \circ \sigma \circ \calC_1 \circ \calE, \label{eq:hybridNetwork}
\end{align}
where $\calL_k : \bbR^{p_k} \to \bbR^{p_{k+1}}, u \mapsto y = \beta_k + \bW_k u$ are the fully connected layers with weight matrices $\bW_k$ and biases $\beta_k$, and $\calC_k$ are the convolutional layers. The parameters of the HNN are collected in $\theta = (b_1,\bK_1,\ldots b_l,\bK_l,\beta_1,\bW_1,\ldots,\beta_{l_2},\bW_{l_2})$. The embedding operator $\calE$, which has essentially been defined in \eqref{eq:image_embedding}, is a form of zero padding that appends zeros to the image such that the convolution operations can be performed. The flattening operator $\calF$ reshapes the output from the convolutional layers into a finite-dimensional vector for the fully connected layers. The flattening operator is defined as
\begin{align*}
	\calF & : \ell_{2e}^{c_l}(\bbN \times \bbN) \to \bbR^{d_l \cdot d_l \cdot c_l}, (u\lsb i_1,i_2 \rsb ) \mapsto y = \begin{pmatrix}
		u\lsb 1,1 \rsb ^\top & \cdots & u\lsb d_l,1 \rsb ^\top & \cdots & u\lsb 1,d_l \rsb ^\top & \cdots & u\lsb d_l,d_l \rsb ^\top
	\end{pmatrix}^\top .
\end{align*}
Note that these operators allow for the evaluations of the convolutional layers to be performed with finite memory and computation. In this fashion, we obtain a hybrid Neural Network as a well-defined mapping $\mathrm{HNN}_{\theta} : \bbR^{d_1\times d_1 \times c_1} \to \bbR^{p_{l_2}}$, where $d_1 \times d_1 \times c_1$ might be the dimensions of an image represented as a vector and $p_{l_2}$ might be the number of classes among which the image should be classified.

We can now make use of dissipativity theory as introduced in Section \ref{sec:3} to obtain upper estimates $\gamma$ on the Lipschitz constant for this class of Neural Networks. To derive a sufficient criterion for $\gamma$ to be an upper bound on the Lipschitz constant, we decompose $\mathrm{HNN}_\theta = \mathrm{FNN}_{\theta_2} \circ \calS_2 \circ \CNN_{\theta_1} \circ \calS_1$ with parameter $\theta = (\theta_1,\theta_2)$ into $\mathrm{FNN}_{\theta_2}$ collecting all the linear layers and $\CNN_{\theta_1}$ collecting all the convolutional layers. We assume that the component $\mathrm{CNN}_{\theta_1}$ is incrementally dissipative with supply rate $s_{\calC}$ and the component $\mathrm{FNN}_{\theta_2}$ is incrementally dissipative with supply rate $s_{\calL}$ defined by
\begin{align*}
	s_\calC(u,y) = \begin{pmatrix}
		y\\
		u
	\end{pmatrix}^\top
	\begin{pmatrix}
		\bQ_C & \bS_C\\
		\bS_C & \bR_C
	\end{pmatrix}
	\begin{pmatrix}
		y\\
		u
	\end{pmatrix}, \qquad
	s_{\calL}(u,y) = 
	\begin{pmatrix}
		y\\
		u
	\end{pmatrix}^\top
	\begin{pmatrix}
		\bQ_L & \bS_L\\
		\bS_L & \bR_L
	\end{pmatrix}
	\begin{pmatrix}
		y\\
		u
	\end{pmatrix}.
\end{align*}
These dissipativity statements can be expressed as
\begin{align}
	\begin{pmatrix}
		y_1^{(l)}\lsb i_1,i_2 \rsb \\
		y_2^{(l)}\lsb i_1,i_2 \rsb 
	\end{pmatrix}
	=
	\begin{pmatrix}
		\CNN_{\theta_1}(u_1^{(1)}\lsb i_1,i_2 \rsb )\\
		\CNN_{\theta_1}(u_2^{(1)}\lsb i_1,i_2 \rsb )
	\end{pmatrix}
	\quad \Rightarrow \quad
	\sum_{i_1,i_2 = 0}^\infty 
	\begin{pmatrix}
		\tilde{y}^{(l)}\lsb i_1,i_2 \rsb \\
		\tilde{u}^{(1)}\lsb i_1,i_2 \rsb
	\end{pmatrix}^\top
	\begin{pmatrix}
		\bQ_C & \bS_C\\
		\bS_C & \bR_C
	\end{pmatrix}
	\begin{pmatrix}
		\tilde{y}^{(l)}\lsb i_1,i_2 \rsb \\
		\tilde{u}^{(1)}\lsb i_1,i_2 \rsb
	\end{pmatrix} \geq 0,  \nonumber
\end{align}
where $\tilde{u}^{(1)}\lsb i_1,i_2 \rsb  = u_1^{(1)}\lsb i_1,i_2 \rsb  - u_2^{(1)}\lsb i_1,i_2 \rsb $, $\tilde{y}^{(l)} =  y_1^{(l)}\lsb i_1,i_2 \rsb  - y_2^{(l)}\lsb i_1,i_2 \rsb $ and
\begin{align}
	\begin{pmatrix}
		y_1^{(l_2)}\\
		y_2^{(l_2)}
	\end{pmatrix}
	=
	\begin{pmatrix}
		\mathrm{FNN}_{\theta_2}(u_1^{(l+1)})\\
		\mathrm{FNN}_{\theta_2}(u_2^{(l+1)})
	\end{pmatrix}
	\quad\Rightarrow\quad
	\begin{pmatrix}
		\tilde{y}^{(l_2)}\\
		\tilde{u}^{(l+1)}
	\end{pmatrix}^\top
	\begin{pmatrix}
		\bQ_L & \bS_L\\
		\bS_L & \bR_L
	\end{pmatrix}
	\begin{pmatrix}
		\tilde{y}^{(l_2)}\\
		\tilde{u}^{(l+1)}
	\end{pmatrix} \geq 0, \nonumber
\end{align}
where $\tilde{y}^{(l_2)} = y_1^{(l_2)} - y_2^{(l_2)}$ and $\tilde{u}^{(l+1)} = \tilde{u}_1^{(l+1)} - \tilde{u}_2^{(l+1)}$. Next, we want to use the incremental dissipativity with respect to $s_{\calL}$ and $s_{\calC}$ to derive a dissipativity statement of the form
\begin{align*}
	\begin{pmatrix}
		y_{1}^{(l_2)}\\
		y_{2}^{(l_2)}
	\end{pmatrix}
	=
	\begin{pmatrix}
		\HNN_{\theta}(u_1^{(1)}\lsb i_1,i_2 \rsb )\\
		\HNN_{\theta}(u_2^{(1)}\lsb i_1,i_2 \rsb )
	\end{pmatrix}
	\quad \Rightarrow \quad
	\gamma^2 \sum_{i_1,i_2 = 0}^{\infty} \| \tilde{u}^{(1)}\lsb i_1,i_2 \rsb \|^2 - \|\tilde{y}^{(l_2)}\|^2 \geq 0,
\end{align*}
which implies an upper bound of $\gamma$ on the Lipschitz constant of $\HNN_{\theta}$. To this end, we further specialize our choice of supply rates $s_\calC$ and $s_\calL$ by requiring $\bR_C = \gamma^2 I$, $\bS_C = 0$, $\bQ_C \preceq 0$, $\bR_L = - \diag (\bQ_C , \ldots , \bQ_C)$, $\bS_L = 0$ and $\bQ_L = - I$. The reason for these choices is simply that it enables us to derive the Lipschitz bound we are seeking from the computation
\begin{align}
    0 &\leq s_{\calL}(\tilde{u}^{(1)},\tilde{y}^{(l_2)}) + \sum_{i_1,i_2 = 0}^{\infty} s_{\calC} (\tilde{u}^{(1)}\lsb i_1,i_2 \rsb,\tilde{y}^{(l)}\lsb i_1,i_2 \rsb ) \nonumber\\
	&\leq \gamma^2 \underbrace{\sum_{i_1,i_2 = 0}^{\infty} \| \tilde{u}^{(1)}\lsb i_1,i_2 \rsb\|^2}_{=\|(\tilde{u}^{(1)}\lsb i_1,i_2 \rsb )\|_{\ell_2}^2}
	+ \underbrace{\sum_{i_1,i_2 = 0}^{\infty} \tilde{y}^{(l)}\lsb i_1,i_2 \rsb ^\top \bQ_C \tilde{y}^{(l)}\lsb i_1,i_2 \rsb 
	+
	(\tilde{u}^{(1)})^\top \bR_L \tilde{u}^{(1)}}_{\overset{(\star)}{\leq}0} - \|\tilde{y}^{(l_2)}\|^2 \nonumber\\
	&\leq
	\gamma^2 \|(\tilde{u}^{(1)}\lsb i_1,i_2 \rsb )\|_{\ell_2}^2 - \|\tilde{y}^{(l_2)}\|^2. \label{eq:LipschitzImplication}
\end{align}
The requirement $\bQ_C \preceq 0$, which we use to conclude $(\star)$, is not very limiting since $\bQ_C \not\preceq 0$ essentially means that the outputs of $\CNN_{\theta_1}$ can grow unbounded even for zero input.
Note that the input time shift by $r$ in $\CNN_{\theta_1}$ does not cause any issue, since it can be eliminated in \eqref{eq:LipschitzImplication} by a change of indices.

Summarizing these findings provides us with a strategy for finding an upper bound on the Lipschitz constant of $\HNN_{\theta}$. We combine the LMI constraint from Theorem \ref{thm:6} for the dissipativity of $\CNN_{\theta_1}$ w.r.t. $s_\calC$ with the LMI constraint from \cite{pauli2021training} for the FNN component and minimize $\gamma^2$ subject to both LMI constraints with free optimization variable $\bQ_C$. This result is formulated in the following theorem.

\begin{theorem}[Lipschitz estimation for CNNs with fully connected layers]
	\label{thm:LipschitzHybridNN}
	Consider a Neural Network of the form \eqref{eq:hybridNetwork}. Let \eqref{eq:error_2Dsys}	be a realization of the incremental dynamics of the convolutional part of the network and let $\bW_1, \ldots, \bW_{l_2}$ be the weight matrices of the linear layers of the linear part of the Neural Network. If there is a solution $(\bQ_C,\gamma^2,\bP_1,\bP_2, \Lambda_C,\Lambda_1,\ldots,\Lambda_{l_2-1})$ to the linear matrix inequalities \eqref{eq:robustDissipativityLMI} and
	\begin{align}
		\left(\begin{array}{cccccc}
			\bR_L            &  -\bW_1^\top \Lambda_1    &   0                 &   \cdots          &         &       \\
			-\Lambda_1 \bW_1    &  2\Lambda_1             & -\bW_2^\top \Lambda_2 &   0         &  \cdots   &  \phantom{\ddots}           \\
			0                 & -\Lambda_2 \bW_2          & 2\Lambda_2          & \ddots            &  \phantom{\ddots} &\vdots       \\
			\vdots            &  0                      & \ddots              & \ddots            &  -\bW_{l_2-1}^\top \Lambda_{l_2-1}  & 0             \\
			&   \vdots                 & \phantom{\ddots}    &  -\Lambda_{l_2-1} \bW_{l_2-1}     &  2\Lambda_{l_2-1}         & -\bW_{l_2}^\top     \\
			&       \phantom{\ddots}   & \cdots              &  0                &  -\bW_{l_2}        & L^2 \bI 
		\end{array}\right) \succeq 0, \label{eq:hybridLMI2}
	\end{align}
	where $\Lambda_C, \Lambda_1,\ldots, \Lambda_{l_2-1}$ are diagonal matrices with non-negative entries, $\bP_1, \bP_2$ are positive semi-definite, and $(\bQ,\bS,\bR)$, $(\bQ^w,\bS^w,\bR^w)$ and $\bR_L$ are parametrized as
    \begin{align*}
        \begin{pmatrix}
            \bQ & \bS\\
            \bS^\top & \bR
        \end{pmatrix}
        &= \begin{pmatrix}
			\bQ_C & 0\\
			0 & \gamma^2 \bI
		\end{pmatrix},
        &
        \begin{pmatrix}
            \bQ^w & \bS^w\\
            (\bS^w)^\top & \bR^w
        \end{pmatrix}
        &=
        \begin{pmatrix}
            -2\Lambda_C & \Lambda_C\\
            \Lambda_C & 0
        \end{pmatrix},
        &
	\bR_L &= 
	\begin{pmatrix}
		-\bQ_C & & \\
		 & \ddots & \\
		 & & -\bQ_C
	\end{pmatrix},
    \end{align*}
    then $\gamma$ is an upper bound on the Lipschitz constant of the mapping \eqref{eq:hybridNetwork}.
\end{theorem}
\begin{proof}
	Consider a Neural Network $\HNN_{\theta} = \FNN_{\theta_2} \circ \calS_2 \circ \CNN_{\theta_1} \circ \calS_1$. According to Theorem \ref{thm:6}, \eqref{eq:robustDissipativityLMI} implies that $\CNN_{\theta_1}$ is incrementally dissipative with respect to the supply rate $s_\calC (u\lsb i_1 ,i_2 \rsb ,y\lsb i_1,i_2 \rsb )$. 
	Furthermore, it is a result of \cite{pauli2021training}, that \eqref{eq:hybridLMI2} implies that $\FNN_{\theta_2}$ is incrementally dissipative with respect to the supply rate $s_{\calL}(u,y) = u^\top \bR_L u - \|y\|^2$. Note that $-\bQ_C$ appears in the left upper diagonal block of \eqref{eq:hybridLMI2}, which implies that $\bQ_C$ is negative semi-definite. 
	Now, we can refer to \eqref{eq:LipschitzImplication}, where we derived that the Lipschitz constant of $\HNN_{\theta}$ is bounded by $\gamma$ if $\CNN_{\theta_1}$ and $\FNN_{\theta_2}$ are dissipative w.r.t. $s_\calC$ and $s_\calL$ and if $\bQ_C \preceq 0$.
\end{proof}

With Theorem \ref{thm:LipschitzHybridNN}, we have an LMI condition for an upper bound on the Lipschitz constant of a Neural Network consisting of linear and convolutional layers. This essentially solves Problem \ref{problem2}. We highlight that this matrix inequality is also linear in the parameter $\gamma^2$, such that we can efficiently optimize over our upper bound.
\section{Lipschitz constant estimation examples for CNNs}
\label{sec:5}

In Sections 2 to 4, we developed methods to represent CNNs as 2-D systems and (semi-definite programs) to estimate their Lipschitz constant. We pointed out the superior scalability of our methodology (for CNNs) in the introduction and we will put it to the test in this section. This improved scalability can already be appreciated for a convolutional Neural Network consisting only of a single convolutional layer.

\begin{example}
	\label{ex:2}
    In this example, we study Lipschitz constant estimation on 100 instances of convolutional layers \eqref{eq:convLayer} defined by randomly sampled $3 \times 3$ kernels $\bK$ with one input and one output channel. All entries of $\bK$ are sampled from independent standard normal distributions. For the estimation of the Lipschitz constant, we consider the following list of methods:
	\begin{enumerate}
		\item[(M1)] As a first possibility, we generate realizations of convolutional layers by forming their transfer functions and computing linear fractional representations using the tool from \cite{biannic2016generalized}. The Lipschitz constant is then estimated by solving
		\begin{align*}
			\minimize_{\bP_1\succeq 0,\bP_2\succeq 0,\gamma^2} ~&~ \gamma^2 \qquad
			\mathrm{s.t.} ~~ \eqref{eq:2DsysLMI},
		\end{align*}
		which is derived from Theorem \ref{thm:5} with the supply rate \eqref{eq:LipschitzSupply}.
            \item[(M2)] Instead of the efficient (in terms of having few states) realization using \cite{biannic2016generalized}, we consider the realization from Example \ref{ex:realization} and solve the same optimization problem as in 1). As we have stressed, this realization is not controllable in the sense of 2-D systems. Since only reachable states will occur in the trajectories of this system, we left and right multiply the matrix inequality \eqref{eq:2DsysLMI} by a basis matrix of the reachable subspace, to consider only reachable states in our dissipativity analysis. We refer to \cite{klamka2005controllability} for the construction of the reachable subspace.
		\item[(M3)] As a third alternative, we can consider the realization as a Fornasini Marchesini system instead of a Roesser system. In this case, we utilize the bounded real lemma for Fornasini Marchesini systems from \cite{agarwal2018bounded} to estimate the Lipschitz constant.
            \item[(M4)] In \cite{fazlyab2019efficient}, a toolbox for the estimation of Lipschitz constants of Neural Networks is provided which we refer to as \emph{LipSDP}. As fourth algorithm for Lipschitz constant estimation, we use this toolbox.
		\item[(M5)] Finally, we consider the operator norm of the Toeplitz operator of \eqref{eq:convLayer} as an approximation (lower bound) on the Lipschitz constant of the convolutional layer.
	\end{enumerate}
	Now we can compare both the computation times and the achieved bounds on the Lipschitz constant of the five methods. These are plotted over the size $d_1$ of the input image in Figure \ref{fig:comparison}.
	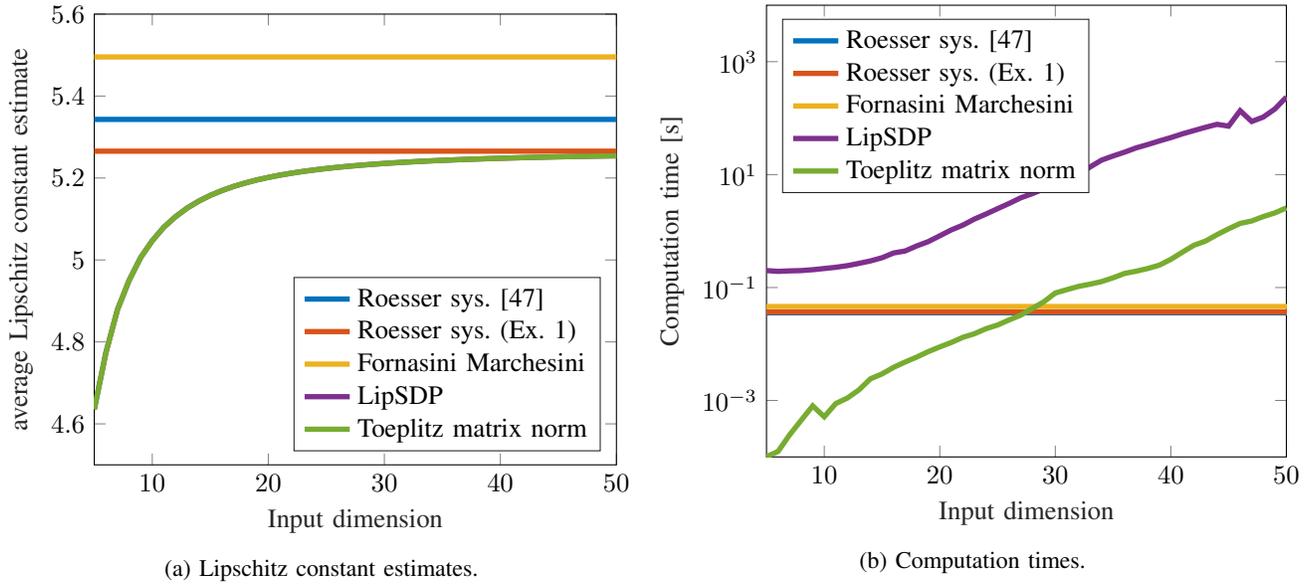
\begin{figure}
		\centering
		\begin{subfigure}{0.48\linewidth}
%
%
\definecolor{mycolor1}{rgb}{0.00000,0.44700,0.74100}%
\definecolor{mycolor2}{rgb}{0.85000,0.32500,0.09800}%
\definecolor{mycolor3}{rgb}{0.92900,0.69400,0.12500}%
\definecolor{mycolor4}{rgb}{0.49400,0.18400,0.55600}%
\definecolor{mycolor5}{rgb}{0.46600,0.67400,0.18800}%
\begin{tikzpicture}

\begin{axis}[%
width=2.733in,
height=2.361in,
at={(2.773in,1.331in)},
scale only axis,
xmin=5,
xmax=50,
xlabel style={font=\color{white!15!black}},
xlabel={Input dimension},
ymin=4.5,
ymax=5.6,
ylabel style={font=\color{white!15!black}},
ylabel={average Lipschitz constant estimate},
axis background/.style={fill=white},
legend style={legend cell align=left, align=left, draw=white!15!black},
legend pos={south east}
]
\addplot [color=mycolor1, line width=2.0pt]
  table[row sep=crcr]{%
5	5.34301903166459\\
50	5.34301903166459\\
};
\addlegendentry{Roesser sys. \cite{biannic2016generalized}}

\addplot [color=mycolor2, line width=2.0pt]
  table[row sep=crcr]{%
5	5.2654626150479\\
50	5.2654626150479\\
};
\addlegendentry{Roesser sys. (Ex. \ref{ex:realization})}

\addplot [color=mycolor3, line width=2.0pt]
  table[row sep=crcr]{%
5	5.49545523135082\\
50	5.49545523135082\\
};
\addlegendentry{Fornasini Marchesini}

\addplot [color=mycolor4, line width=2.0pt]
  table[row sep=crcr]{%
5	4.63548527401232\\
6	4.77334542366022\\
7	4.8786592745162\\
8	4.94991776697037\\
9	5.00596432996976\\
10	5.0469600656478\\
11	5.07997897486533\\
12	5.10546974054451\\
13	5.12648102995916\\
14	5.14338269311224\\
15	5.15756866411027\\
16	5.16928394201349\\
17	5.17933884216326\\
18	5.18776230686466\\
19	5.19516863531493\\
20	5.20143246950189\\
21	5.2070024408604\\
22	5.21179314047574\\
23	5.21608243506045\\
24	5.21983049806355\\
25	5.22320082915659\\
26	5.22619312966348\\
27	5.22888961189753\\
28	5.23131278905955\\
29	5.2335052382199\\
30	5.2354920478978\\
31	5.23730641831763\\
32	5.2389497648315\\
33	5.24046793574991\\
34	5.24184178350134\\
35	5.24312515644878\\
36	5.2442855910288\\
37	5.24537905858716\\
38	5.24636837794636\\
39	5.24730848234933\\
40	5.24816111171027\\
41	5.24897221803246\\
42	5.24971339446855\\
43	5.25041712100949\\
44	5.25106557422743\\
45	5.25167974650554\\
46	5.25225055610633\\
47	5.25278998103313\\
48	5.25329497622898\\
49	5.25377123481122\\
50	5.25421990983372\\
};
\addlegendentry{LipSDP}

\addplot [color=mycolor5, line width=2.0pt]
  table[row sep=crcr]{%
5	4.63548527401331\\
6	4.77334542372555\\
7	4.87865927460164\\
8	4.94991776701399\\
9	5.00596433001634\\
10	5.046960065654\\
11	5.07997897487823\\
12	5.1054697405451\\
13	5.12648102981906\\
14	5.1433826929634\\
15	5.15756866410376\\
16	5.1692839417934\\
17	5.17933884210431\\
18	5.18776230658443\\
19	5.19516863499731\\
20	5.20143246898184\\
21	5.20700244043689\\
22	5.21179314001567\\
23	5.21608243452359\\
24	5.21983049744457\\
25	5.22320082843708\\
26	5.22619312885803\\
27	5.22888961095318\\
28	5.2313127884252\\
29	5.23350523763466\\
30	5.23549204709055\\
31	5.23730641765893\\
32	5.23894976412093\\
33	5.24046793511472\\
34	5.24184178311309\\
35	5.24312515605839\\
36	5.24428558858149\\
37	5.24537905838974\\
38	5.24636837774298\\
39	5.24730848193944\\
40	5.2481611111753\\
41	5.24897221766379\\
42	5.24971339413338\\
43	5.25041711947104\\
44	5.25106557391463\\
45	5.25167974618871\\
46	5.25225055578602\\
47	5.25278998076145\\
48	5.25329497582505\\
49	5.25377123450623\\
50	5.25421990937512\\
};
\addlegendentry{Toeplitz matrix norm}

\end{axis}
\end{tikzpicture}%
			\caption{Lipschitz constant estimates.}
		\label{fig:comparisona}
		\end{subfigure}
		\begin{subfigure}{0.48\linewidth}
%
%
\definecolor{mycolor1}{rgb}{0.00000,0.44700,0.74100}%
\definecolor{mycolor2}{rgb}{0.85000,0.32500,0.09800}%
\definecolor{mycolor3}{rgb}{0.92900,0.69400,0.12500}%
\definecolor{mycolor4}{rgb}{0.49400,0.18400,0.55600}%
\definecolor{mycolor5}{rgb}{0.46600,0.67400,0.18800}%
\begin{tikzpicture}

\begin{axis}[%
width=2.721in,
height=2.366in,
at={(0.758in,0.481in)},
scale only axis,
xmin=5,
xmax=50,
xlabel style={font=\color{white!15!black}},
xlabel={Input dimension},
ymode=log,
ymin=0.0001,
ymax=10000,
yminorticks=true,
ylabel style={font=\color{white!15!black}},
ylabel={Computation time [s]},
axis background/.style={fill=white},
legend style={legend cell align=left, align=left, draw=white!15!black},
legend pos={north west}
]
\addplot [color=mycolor1, line width=2.0pt]
  table[row sep=crcr]{%
5	0.035683987\\
100	0.035683987\\
};
\addlegendentry{Roesser sys. \cite{biannic2016generalized}}

\addplot [color=mycolor2, line width=2.0pt]
  table[row sep=crcr]{%
5	0.036864945\\
100	0.036864945\\
};
\addlegendentry{Roesser sys. (Ex. \ref{ex:realization})}

\addplot [color=mycolor3, line width=2.0pt]
  table[row sep=crcr]{%
5	0.046151592\\
100	0.046151592\\
};
\addlegendentry{Fornasini Marchesini}

\addplot [color=mycolor4, line width=2.0pt]
  table[row sep=crcr]{%
5	0.199934648\\
6	0.193362088\\
7	0.197237231\\
8	0.199766785\\
9	0.206922236\\
10	0.217282508\\
11	0.228366988\\
12	0.244169687\\
13	0.267801603\\
14	0.295941291\\
15	0.33704117\\
16	0.409403614\\
17	0.443929132\\
18	0.545607886\\
19	0.650356222\\
20	0.824895686\\
21	1.049421538\\
22	1.274042278\\
23	1.632732009\\
24	2.002880154\\
25	2.493147817\\
26	3.082450901\\
27	3.877335269\\
28	4.59071467\\
29	5.588206046\\
30	6.558964928\\
31	7.680703325\\
32	10.092334606\\
33	13.615993323\\
34	18.040608523\\
35	21.38604418\\
36	25.010727636\\
37	29.750695127\\
38	34.082257996\\
39	39.418329368\\
40	45.24654839\\
41	52.656143285\\
42	60.258559491\\
43	68.843219387\\
44	77.774321232\\
45	72.182588825\\
46	135.509503492\\
47	87.375243547\\
48	104.441980251\\
49	145.017673821\\
50	235.961981572\\
};
\addlegendentry{LipSDP}

\addplot [color=mycolor5, line width=2.0pt]
  table[row sep=crcr]{%
5	0.000100877\\
6	0.000124975\\
7	0.000245806\\
8	0.000444735\\
9	0.00080273\\
10	0.000513919\\
11	0.000882084\\
12	0.001103329\\
13	0.001536932\\
14	0.002420528\\
15	0.002953851\\
16	0.003859414\\
17	0.004794844\\
18	0.005859359\\
19	0.007322054\\
20	0.008879602\\
21	0.010645562\\
22	0.013119887\\
23	0.015203146\\
24	0.018636649\\
25	0.021672808\\
26	0.026655241\\
27	0.032385417\\
28	0.042489161\\
29	0.055746293\\
30	0.08002152\\
31	0.091755404\\
32	0.104639913\\
33	0.114712737\\
34	0.127377163\\
35	0.148776471\\
36	0.17786732\\
37	0.194623135\\
38	0.217360505\\
39	0.248992608\\
40	0.315627891\\
41	0.423012929\\
42	0.561102151\\
43	0.664942204\\
44	0.864310996\\
45	1.096402159\\
46	1.365885529\\
47	1.510244622\\
48	1.806421525\\
49	2.096537518\\
50	2.572288371\\
};
\addlegendentry{Toeplitz matrix norm}

\addplot [color=mycolor4, dashed, line width=2.0pt, forget plot]
  table[row sep=crcr]{%
51	294.798923504741\\
52	366.718041315372\\
53	454.289514812562\\
54	560.524657546389\\
55	688.941115149093\\
56	843.636206440625\\
57	1029.36924492135\\
58	1251.65374668504\\
59	1516.86050338696\\
60	1832.3325758691\\
61	2206.51334555526\\
62	2649.08884696534\\
63	3171.14569584695\\
64	3785.34602367345\\
65	4506.12093080508\\
66	5349.88407765205\\
67	6335.26714591807\\
68	7483.37902064543\\
69	8818.09066853941\\
70	10366.3478191351\\
71	12158.5136930014\\
72	14228.7441655792\\
73	16615.3979066493\\
74	19361.4841940508\\
75	22515.1512663625\\
76	26130.2182530467\\
77	30266.753902298\\
78	34991.705516765\\
79	40379.5817056936\\
80	46513.1927691166\\
81	53484.4527457576\\
82	61395.2473815834\\
83	70358.3725107043\\
84	80498.5475848574\\
85	91953.5093422864\\
86	104875.190871745\\
87	119430.991602874\\
88	135805.144040638\\
89	154200.183359135\\
90	174838.526279235\\
91	197964.165975428\\
92	223844.490090331\\
93	252772.229280753\\
94	285067.544077462\\
95	321080.258212055\\
96	361192.24694901\\
97	405819.989359389\\
98	455417.293885104\\
99	510478.20696948\\
100	571540.114971437\\
};
\addplot [color=mycolor5, dashed, line width=2.0pt, forget plot]
  table[row sep=crcr]{%
51	3.01944629959861\\
52	3.53332246931629\\
53	4.12229480720996\\
54	4.79560314103485\\
55	5.56341907179161\\
56	6.43692004509805\\
57	7.42836778893891\\
58	8.55119128877541\\
59	9.82007447443074\\
60	11.2510487966077\\
61	12.8615908743402\\
62	14.6707253981325\\
63	16.6991334769959\\
64	18.9692666210565\\
65	21.5054665548734\\
66	24.3340910600825\\
67	27.4836460494558\\
68	30.9849240779546\\
69	34.8711494998361\\
70	39.1781304843758\\
71	43.9444181062568\\
72	49.2114727301883\\
73	55.0238379128194\\
74	61.429322048527\\
75	68.4791879891775\\
76	76.2283508714802\\
77	84.7355843890788\\
78	94.0637357500597\\
79	104.279949564089\\
80	115.455900906929\\
81	127.668037813641\\
82	140.997833455303\\
83	155.532048257662\\
84	171.363002223663\\
85	188.58885772537\\
86	207.31391303438\\
87	227.648906863355\\
88	249.711334194921\\
89	273.625773677729\\
90	299.524226873048\\
91	327.546469638886\\
92	357.840415942174\\
93	390.562494393174\\
94	425.87803779986\\
95	463.961686043614\\
96	504.997802581194\\
97	549.180904881529\\
98	596.716109109516\\
99	647.819589372607\\
100	702.719051849575\\
};
\end{axis}
\end{tikzpicture}%
			\caption{Computation times.}
		\label{fig:comparisonb}
		\end{subfigure}
		\caption{This figure presents the average Lipschitz constant estimates and the average computation times of five methods for estimating the Lipschitz constant of a convolutional layer. These averages are calculated over a set of 100 randomly generated convolutional layers. Note that, in Figure \ref{fig:comparisona}, the plot for the norm of the Toeplitz operator occludes the plot of LipSDP (both equal the exact Lipschitz constant of the convolutional layer).}
		\label{fig:comparison}
	\end{figure}
    As we can see in the figure, the methods M1, M2, and M3 provide uniform bounds over the input size $d_1$ and also enjoy uniform computation time statistics. That is because these approaches are based on 2-D system realizations and their Lipschitz constant bounds apply for any size of the input. The methods M4 and M5, on the other hand, depend on the size of the input and, consequently, their computation time grows for larger $d_1$. It should also be stressed that for a fixed given $d_1$, these methods provide the more accurate Lipschitz constant bound.
    Of course, the case of one single layer is an extreme and artificial experiment, which favors a 2-D system approach and is unfavorable for LipSDP. On the other hand, the computational costs for estimating the Lipschitz constant of a convolutional Neural Network for images of size $d_1$ is definitely going to be lower bounded by the numbers provided in Figure \ref{fig:comparison}. Hence, we can see that converting convolutional layers into linear layers becomes quite unfavorable already for input sizes above $d_1\times d_1 = 50 \times 50$, which are still quite moderate resolutions for images. (Empirical estimations indicated a growth of computational effort of the order $d_1^{11}$ for LipSDP and $d_1^8$ for the Toeplitz norm approach.)
    Consequently, we draw the conclusion that methods tailored to convolutional layers are required for the dissipativity-based Lipschitz constant estimation of Neural Networks.
    
    Finally, we would like to point towards an interesting phenomenon, which we can observe in Figure \ref{fig:comparison}, namely, that different 2-D system realizations of the same transfer functions can lead to different bounds on the Lipschitz constant. Empirically, the Fornasini Marchesini realization leads to more conservative estimates than the Roesser realization (which is one reason why we decided for Roesser systems), and the realization from Example \ref{ex:realization}, with more redundancy, results in more accurate estimates than the smaller Roesser realization. Indeed, the Lipschitz constant bounds produced by the redundant realization appear to be exact. This exactness is an interesting observation that should be subject to further research. What we generally observed was that lifting the 2-D system to realizations with additional states leads to improved estimates on the $\ell_2$-gain. A result, which is also reported in \cite{scherer2016lossless}, where an asymptotically exact relaxation hierarchy for $H_\infty$-synthesis of 2-D systems with a continuous and a discrete axis is proposed.
\end{example}

Example \ref{ex:2}, where we consider a CNN consisting of only one convolutional layer, shows that our 2-D system approach can be advantageous for large input images. In the next example, we explore whether this approach is also competitive for CNNs with multiple layers trained for image classification.

\begin{figure}
	\centering
	\includegraphics[width = 0.75\linewidth]{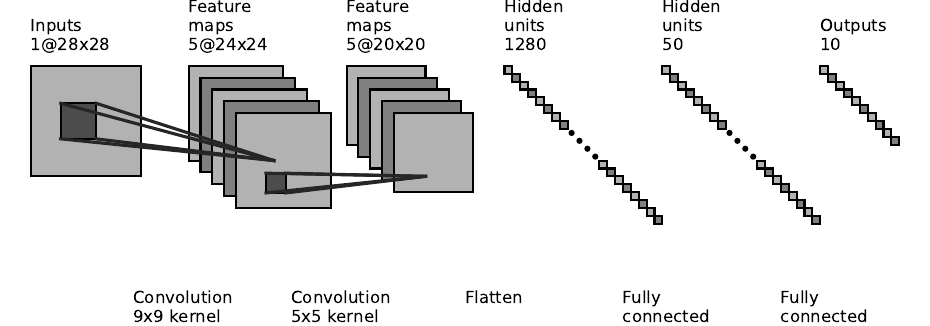}
	\caption{This figure shows the Neural Network architecture we trained on the MNIST data set. The input image size is $28\times 28$. We then apply one $9\times 9$ convolution and one $5\times 5$ convolution with five channels each, which maps the input image to an image of the size $16\times 16$ with five channels. The output of the second convolutional layer is then flattened and processed by two fully connected layers with 50 neurons and 10 output neurons. Between each of these layers there is one layer of ReLU activation functions.}
	\label{fig:architecture}
\end{figure}

\begin{example}[Lipschitz constant estimation for CNNs on MNIST]
	\label{ex:3}
	In this example, we study the Lipschitz constant estimation for a Neural Network with convolutional and fully connected layers trained to classify the hand-written digits of the MNIST dataset. The Neural Network structure we consider consists of two convolutional layers with $9\times 9$ and $5 \times 5$ convolution kernels and two linear layers with $50$ neurons and $10$ neurons. This structure is depicted in Figure \ref{fig:architecture}. Training this network on the MNIST dataset for 14 epochs with stochastic gradient descent and a learning rate of 0.01 produces a classification error of $0.88\%$ on the test set (the best Neural Network achieves an error of $0.16\%$ \cite{mazzia2021efficient}). For the trained Neural Network, we carried out the Lipschitz constant estimation using Theorem \ref{thm:LipschitzHybridNN} and using the method \emph{LipSDP} introduced in \cite{fazlyab2019efficient}. The required 2-D system representations have been generated as linear fractional representations of the transfer functions of the convolution kernel using the toolbox developed in \cite{biannic2016generalized}. With this 2-D system representation, we solve the optimization problem
	\begin{align}
		\minimize_{\bP_1\succeq 0,\bP_2\succeq 0,\Lambda_C \geq 0,\Lambda_1,\ldots,\Lambda_{l_2},\gamma^2} ~&~ \gamma^2 \label{eq:OptimizeHybrid}\\
		\mathrm{s.t.} ~&~ \eqref{eq:robustDissipativityLMI}, \eqref{eq:hybridLMI2} \nonumber
	\end{align}
	which is a convex optimization problem. Due to Theorem \ref{thm:LipschitzHybridNN}, this yields an upper bound on the Lipschitz constant of the Neural Network under consideration, which we feature in Table \ref{tab:hybridLipschitzConstant} together with the computation time.
	\begin{table}
		\centering
		\begin{tabular}{@{}llll@{}}
			\toprule
			Solver & Input size & Solver time & Lipschitz constant upper bound \\
			\midrule
			2-D system approach & $d_1 = 14$ & 82.7 s & 23.9 \\
			LipSDP & $d_1 = 14$ & 334 s & 12.0 \\
			Toeplitz norm & $d_1 = 14$ & 0.0108 s & 43.0 \\[2mm]
			2-D system approach & $d_1 = 28$ & 1550 s & 18.9 \\
			LipSDP & $d_1 = 28$ & out of memory & - \\
			Toeplitz norm & $d_1 = 28$ & 0.591 s & 47.3 \\
			\bottomrule
		\end{tabular}
		\caption{Lipschitz constant upper bound and solver time in seconds for solving the optimization problem \eqref{eq:OptimizeHybrid} using Mosek.}
		\label{tab:hybridLipschitzConstant}
	\end{table}
    Because LipSDP is not able to provide a Lipschitz constant estimate for this example due to excessive storage requirements, we also consider downsampled images with $d_1 \times d_1 = 14 \times 14$ pixels. (Note that this decreases the number of neurons from 5434 to 686.) Still, for this reduced example, the 2-D system approach is significantly faster than LipSDP. On the other hand, we see in Table \ref{tab:hybridLipschitzConstant}, that the Lipschitz constant estimates of LipSDP are tighter than those provided by \eqref{eq:robustDissipativityLMI}. We assume this to be the case, because (particularly in the $14\times 14$ images), the size of the convolution kernel is not small compared to the size of the image, and because the transition from convolutional to linear layers might introduce some conservatism.
\end{example}

Example \ref{ex:3} considers a setting that is much closer to practical applications than Example \ref{ex:2}, namely digit classification. The results show that already for the $28 \times 28$ images of the MNIST dataset, established methods (LipSDP) for Lipschitz constant estimation based on robust control get into trouble. The proposed tools based on 2-D systems prove superior computational feasibility also for this case of a multi-layer CNN including linear layers. Note, that in Example \ref{ex:3} a ReLU Neural Network is considered for solving the MNIST classification problem. While being standard in Deep Learning, the ReLU activation function is not considered competitive for learning Lipschitz continuous functions \cite{bethune2022pay}, where \emph{group sort activation functions} \cite{chernodub2016norm,anil2019sorting,tanielian2021approximating} should be more appropriate.
\section{Conclusion}

In this exposition, we show how convolutional Neural Networks can be represented as 2-D Lur´e systems, i.e., as a hybrid of Lur´e systems and 2-D systems from control. We believe that this representation is quite valuable, on the one hand, since CNNs are by now among the most important applications for two-dimensional filters (and two-dimensional filters are a core interest in 2-D systems research), and, on the other hand, since the research stream of estimating the Lipschitz constant of Neural Networks using robust control has been lacking an efficient representation of CNNs, which can be provided by our 2-D systems representation. As further results of our research, we apply dissipativity analysis to 2-D Lur´e systems and derive SDPs for Lipschitz constant estimation of CNNs. Our examples show that the efficient representation of CNNs enables this Lipschitz constant estimation for CNNs with a considerably larger amount of neurons than previously possible (in terms of computational tractability).





\bibliographystyle{IEEEtran}
\bibliography{references}







\end{document}